\documentclass[12pt,francais]{smfart}

\usepackage[dvips]{epsfig}

\usepackage[francais]{babel}

\usepackage{smfthm}

\newcommand{\C}{\mathbb C}
\renewcommand{\P}{\mathbb P}
\newcommand{\N}{\mathbb N}

\renewcommand{\L}{\mathbb L}

\newcommand{\G}{\mathbb G}

\newtheorem{theorem}{Th\'eor\`eme}[section]
\newtheorem{lemma}[theorem]{Lemme}
\newtheorem{proposition}[theorem]{Proposition}
\newtheorem{corollary}[theorem]{Corollaire}
\newtheorem{rappel}[theorem]{}
\newtheorem{problem}[theorem]{Probl\`eme}
\newtheorem{propriete}[theorem]{Propri\'et\'e}

\theoremstyle{definition}
\newtheorem{definition}[theorem]{D\'efinition}
\newtheorem{remark}[theorem]{Remarque}

\newtheorem{notation}[theorem]{Notation}
\newtheorem{question}[theorem]{Question}
\newtheorem{example}[theorem]{Exemple}
\newtheorem{conjecture}[theorem]{Conjecture}
\newtheorem{exercice}[theorem]{Exercice}

\newcommand{\cal}{\mathcal}
\newcommand{\wh}{\widehat}
\newcommand{\lan}{\langle}
\newcommand{\ran}{\rangle}

\newcommand{\bi}{\begin{itemize}}
\newcommand{\ei}{\end{itemize}}
\newcommand{\be}{\begin{enumerate}}
\newcommand{\ee}{\end{enumerate}}

\newcommand{\bpf}{\begin{proof}}
\newcommand{\epf}{\end{proof}}

\newcommand{\bpro}{\begin{propriete}}
\newcommand{\epro}{\end{propriete}}

\newcommand{\bt}{\begin{theorem}}
\newcommand{\et}{\end{theorem}}
\newcommand{\brap}{\begin{rappel}}
\newcommand{\erap}{\end{rappel}}
\newcommand{\bnt}{\begin{notation}}
\newcommand{\ent}{\end{notation}}
\newcommand{\bd}{\begin{definition}}
\newcommand{\ed}{\end{definition}}
\newcommand{\ble}{\begin{lemma}}
\newcommand{\ele}{\end{lemma}}
\newcommand{\bpr}{\begin{proposition}}
\newcommand{\epr}{\end{proposition}}
\newcommand{\bre}{\begin{remark}}
\newcommand{\ere}{\end{remark}}
\newcommand{\bco}{\begin{corollary}}
\newcommand{\eco}{\end{corollary}}
\newcommand{\beq}{\begin{equation}}
\newcommand{\eeq}{\end{equation}}
\newcommand{\bq}{\begin{question}}
\newcommand{\eq}{\end{question}}
\newcommand{\bp}{\begin{problem}}
\newcommand{\ep}{\end{problem}}
\newcommand{\beqn}{\begin{eqnarray*}}
\newcommand{\eeqn}{\end{eqnarray*}}
\newcommand{\bex}{\begin{example}}
\newcommand{\eex}{\end{example}}
\newcommand{\ber}{\begin{exercice}}
\newcommand{\eer}{\end{exercice}}

\newcommand{\sct}{\section}
\newcommand{\ssct}{\subsection}

\newcommand{\sk}{\smallskip}
\newcommand{\bk}{\bigskip}
\newcommand{\nk}{\noindent}
\newcommand{\pl}{\partial}
\newcommand{\ov}{\overline}
\newcommand{\fr}{\frac}
\newcommand{\bcj}{\begin{conjecture}}
\newcommand{\ecj}{\end{conjecture}}

\newcommand{\bck}{\backslash}

\newcommand{\Si}{\Sigma}
\newcommand{\s}{\sigma}

\renewcommand{\t}{\theta}
\renewcommand{\a}{\alpha}
\renewcommand{\b}{\beta}
\newcommand{\g}{\gamma}

\renewcommand{\l}{\lambda}
\renewcommand{\O}{\Omega}
\renewcommand{\k}{\kappa}
\renewcommand{\o}{\omega}

\title[Sur l'alg\'ebrisation des tissus]{Sur l'alg\'ebrisation des \\ tissus de rang maximal}

\author[L. Pirio et J.-M. Tr\'epreau]{Luc Pirio et Jean-Marie Tr\'epreau}

\begin{abstract}
Soit $n\geq 2$, $r\geq 2$ et $d\geq  (r+1)(n-1)+2$ des entiers.
Nous montrons qu'un $d$-tissu de codimension $r$ sur un germe de vari\'et\'e 
de dimension $rn$
est alg\'ebrisable au sens classique s'il est de rang maximal,
sauf peut-\^etre lorsque $n\geq 3$ et $d=(r+2)(n-1)+1$.
Dans ce dernier cas, nous montrons qu'il est alg\'ebrisable, mais 
en un sens g\'en\'eralis\'e.
\end{abstract}

\begin{document}

\maketitle

Cet article concerne ce que nous appellerons la th\'eorie ab\'elienne 
des tissus. 
Nous \'etendons aux
tissus de codimension quelconque le th\'eor\`eme
d'alg\'ebrisation des tissus de codimension un et de rang 
maximal de Tr\'epreau \cite{Tr}. 
Nous r\'epondons ainsi au probl\`eme pos\'e 
par Chern et Griffiths dans \cite {C-Gr2}, \cite{C-Gr3}.
L'objet de cette pr\'esentation g\'en\'erale est 
d'\'enoncer notre r\'esultat, apr\`es les 
quelques pr\'eliminaires n\'ecessaires \`a sa compr\'ehension.
Les notions introduites  rapidement ci-dessous 
sont d\'efinies pr\'ecis\'ement dans les deux premi\`eres 
sections du Chapitre~1.

\sk
Un $d$-tissu $\cal{T}=\{\cal{F}_1,\ldots,\cal{F}_d\}$ de codimension 
$r$ sur un germe $(M,x_0)$ de vari\'et\'e analytique 
complexe lisse est une famille de $d$ feuilletages r\'eguliers
de codimension $r$  de $(M,x_0)$, v\'erifiant une condition 
de \og position g\'en\'erale \fg.

L'exemple suivant est fondamental.
Soit $Z\subset \P^{r+n-1}$ une vari\'et\'e projective  
de dimension $r$ et de degr\'e $d$, telle qu'un $(n-1)$-plan
g\'en\'eral rencontre $Z$ en $d$ points qui sont en position
g\'en\'erale dans ce plan. La dualit\'e 
projective permet de lui associer un $d$-tissu 
alg\'ebrique $\cal{T}_Z$ de codimension $r$ sur la grassmannienne $\G_{r,n}$ 
des $(n-1)$-plans de $\P^{r+n-1}$.

Au point g\'en\'erique de la grassmannienne, ce tissu alg\'ebrique 
induit un germe de tissu, au sens du paragraphe pr\'ec\'edent.
Nous dirons que $\cal{T}_Z$, ou un germe de tissu d\'efini par $\cal{T}_Z$,
est un tissu {\em alg\'ebrique 
grassmannien}.

\sk
La grassmannienne $\G_{r,n}$
est une vari\'et\'e de dimension $rn$. Nous supposons 
maintenant que le germe $(M,x_0)$ est de dimension $rn$ avec $n\geq 2$.

\sk
La th\'eorie ab\'elienne des tissus na\^{\i}t en 1932 avec l'introduction 
par Blaschke de la notion de {\em relation ab\'elienne d'un tissu}
dans le cas des tissus plans, notion \'etendue peu apr\`es aux tissus de
codimension un et beaucoup plus tard, dans Griffiths \cite{Gr} et Chern-Griffiths \cite{C-Gr3}, au cas qui nous importe ici.

\sk
Pour $j=1,\ldots,d$, soit $\Omega_j$ une normale d\'efinissante
du feuilletage $\cal{F}_j$. Une relation ab\'elienne du tissu $\cal{T}$ est un 
$d$-uplet de {\em normales ferm\'ees} $(f_1\O_1,\ldots,f_d\O_d)$,  {\em i.e.}
les $f_j$ sont des fonctions analytiques sur $(M,x_0)$ telles que $d(f_j \Omega_j)= 0$,
v\'erifiant la relation  $\sum_{j=1}^d f_j \Omega_j = 0$.
{\em Le rang} du tissu $\cal{T}$ est la dimension du $\C$-espace vectoriel
de ses relations ab\'eliennes. 

\sk
Des bornes sur le rang d'un tissu ont \'et\'e donn\'ees par Bol (1932)
pour $(r,n)=(1,2)$, par Chern (1935) pour $r=1$ avec $n$ quelconque. Dans \cite{C-Gr3}
Chern et Griffiths consid\`erent le cas g\'en\'eral et montrent,
et cette borne est optimale,
que le rang du tissu $\cal{T}$ est au plus \'egal au nombre  
$$
\rho_{r,n}(d) = \sum_{h=0}^{+\infty} \, { r-1 +h \choose r-1 } \times \max\,(d-(r+h)(n-1)-1,0).
$$

Revenons au cas d'un tissu alg\'ebrique sur la 
grassmannienne $\G_{r,n}$ d\'efini par une sous-vari\'et\'e projective $Z$ de $\P^{r+n-1}$,
de dimension $r$ et de degr\'e $d$, que nous supposons, pour simplifier,
irr\'eductible et lisse.

Il r\'esulte de la borne pr\'ec\'edente et 
de la version g\'en\'erale du th\'eor\`eme d'Abel
donn\'ee par Griffiths  \cite{Gr},
que le genre g\'eom\'etrique de $Z$ est major\'e par
$\rho_{r,n}(d)$. C'est une g\'en\'eralisation int\'eressante 
de la borne classique de Castelnuovo pour le genre des courbes de $\P^n$
de degr\'e donn\'e.
Peu apr\`es, Harris \cite{Ha} a caract\'eris\'e,
mais essentiellement
dans le cadre des vari\'et\'es irr\'eductibles et lisses,
les vari\'et\'es qui sont de genre g\'eom\'etrique maximal.

\sk
C'est le probl\`eme de la r\'eciproque, pos\'e \`a cette \'epoque,
que nous r\'esolvons ici. Notre r\'esultat principal est le suivant.

\bk
\nk
{\bf {\em Th\'eor\`eme}}.---
{\em Soit $\cal{T}$ un $d$-tissu de codimension $r$ et
de rang maximal $\rho_{r,n}(d)$, avec 
$r\geq 2$ et $d> (r+1)(n-1)+1$.

Si de plus $\;d\neq (r+2)(n-1)+1\;$ ou si $\;n=2$, il est localement
isomorphe \`a un germe de tissu alg\'ebrique grassmannien.

Si $d=(r+2)(n-1)+1$ et $n\geq 3$, il est localement isomorphe 
\`a un germe d'un tissu alg\'ebrique, \'eventuellement non
grassmannien.}

\sk
Si $d\leq r(n-1)+1$, le nombre $\rho_{r,n}(d)$ est nul 
et la question de la r\'eciproque ne se pose pas.
M\^eme si la premi\`ere condition de l'\'enonc\'e 
\'ecarte les cas o\`u $r(n-1)+1 < d \leq (r+1)(n-1)+1$,
qui m\'eriteraient d'\^etre \'etudi\'es, elle est naturelle
pour des raisons g\'eom\'etriques et compte tenu de la 
m\'ethode que nous utilisons.

Si $n=2$ ou si $d \neq (r+2)(n-1)+1$, nous obtenons le r\'esultat attendu.

\sk
Le cas  $d = (r+2)(n-1)+1$ avec $n\geq 3$ est une surprise.
Il faut alors introduire une classe plus g\'en\'erale de 
{\em tissus alg\'ebriques} pour avoir l'\'enonc\'e ci-dessus.
Le probl\`eme de la d\'etermination des $d$-tissus de rang
maximal pour cette valeur de $d$ reste ouvert. Nous donnerons 
dans un article \`a venir \cite{P-Tr2}, pour quelques valeurs de $r$ et de $n$,
des exemples de $d$-tissus de rang maximal qui sont alg\'ebrisables,
mais qui ne sont pas isomorphes \`a des germes
de tissus alg\'ebriques grassmanniens.

\sk
L'\'enonc\'e pr\'ec\'edent est pr\'ecis\'e, mis en 
en perspective et discut\'e dans l'introduction qui suit.

\sct{Introduction}

La th\'eorie ab\'elienne des tissus na\^it avec la d\'efinition par Blaschke
de la notion de relation ab\'elienne d'un tissu plan. Les r\'esultats majeurs 
obtenus entre 1932 et 1935, en particulier  par Blaschke, Howe et Bol,
sont pr\'esent\'es dans la derni\`ere partie du livre {\em Geometrie der Gewebe} 
de Blaschke et Bol \cite{BB}, paru en 1938.

\sk
Les fondateurs de la th\'eorie font deux emprunts essentiels \`a la 
litt\'erature classique. Le premier \`a Poincar\'e \cite{Po} qui, 
cherchant une nouvelle d\'emonstration du th\'eor\`eme de Lie sur 
les surfaces de double translation,  introduit une id\'ee
qui restera  l'id\'ee fondamentale pour r\'esoudre les probl\`emes
d'alg\'ebrisation des tissus de rang maximal. 

Le second \`a Darboux \cite{Da} qui reprend cette id\'ee mais remplace
la seconde partie peu claire de l'analyse de Poincar\'e par un argument tr\`es ing\'enieux,
le prototype de ce que l'on appelle le th\'eor\`eme d'Abel inverse,
obtenant une d\'emonstration g\'eom\'etrique limpide du r\'esultat de Lie.

\sk 
Dans ce chapitre, nous pr\'esentons d'abord le probl\`eme et les fondements de la 
th\'eorie, avec comme r\'ef\'erences principales Chern-Griffiths~\cite{C-Gr3}
et Griffiths~\cite{Gr}.
Nous rappelons ensuite le principe de 
la {\em m\'ethode canonique}. La plus grande partie de cet article est
consacr\'ee \`a la construction, par cette m\'ethode, de 
{\em la vari\'et\'e de Blaschke} d'un tissu de rang maximal
et \`a en \'etablir les propri\'et\'es essentielles. Elles sont 
d\'ecrites dans le Th\'eor\`eme~\ref{Th-X0}.  

\sk
Le probl\`eme de l'alg\'ebrisation est alors r\'eduit 
\`a un probl\`eme de g\'eom\'etrie projective.
Nos r\'esultats principaux, les Th\'eor\`emes  \ref{Th1A}, \ref{Th1B} et \ref{Th1C}
reposent pour une part importante sur la solution 
que nous avons apport\'ee \`a  ce probl\`eme dans Pirio-Tr\'epreau~\cite{P-Tr}.

\sk
Nous terminerons par quelques rappels sur l'histoire
du probl\`eme que nous r\'esolvons ici.

\ssct{Tissus et relations ab\'eliennes}

Soit $r\geq 1$ et $n\geq 2$ des entiers. Un {\em feuilletage de type $(r,n)$}
est un  feuilletage r\'egulier $\cal{F}$,  de codimension~$r$,
sur un germe $(M,x_0)$ de vari\'et\'e analytique (complexe) lisse de dimension~$rn$.
D'une fa\c{c}on g\'en\'erale, les objets consid\'er\'es sur $(M,x_0)$
seront analytiques, sans qu'on le pr\'ecise.

Souvent, en particulier dans les \'enonc\'es principaux, nous supposerons
$r\geq 2$. La raison pour cela est que le probl\`eme \'etudi\'e 
ici est d\'ej\`a r\'esolu pour $r=1$ et qu'il y a 
quelques diff\'erences de traitement entre le cas $r=1$ et le cas $r\geq 2$.

\sk
Le feuilletage $\cal{F}$ est par exemple d\'efini par un syst\`eme diff\'erentiel 
$$
du_a=0, \qquad a=1,\ldots,r,
$$
o\`u les $u_a$ sont des fonctions sur $(M,x_0)$ 
dont les diff\'erentielles en $x_0$ sont 
lin\'eairement ind\'ependantes. Les feuilles du feuilletage $\cal{F}$
sont les vari\'et\'es int\'egrales de  ce syst\`eme. 

\sk
Nous appelons {\em normale de $\cal{F}$} toute
$r$-forme $\phi = f du_1\wedge \cdots \wedge du_r$ , o\`u $f$ est une fonction sur $(M,x_0)$.
La normale $\phi$ est {\em ferm\'ee} si $d\phi=0$.
C'est une {\em normale  g\'en\'eratrice} si $f$ ne s'annule pas.
Ces notions ne d\'ependent pas du syst\`eme de fonctions choisi
pour pr\'esenter $\cal{F}$.

\bk
Nous consid\'erons maintenant des familles de feuilletages.  
Soit $d\geq 1$ un entier et $\cal{T}=\{\cal{F}_1,\ldots,\cal{F}_d\}$
une famille de feuilletages de type $(r,n)$ sur le germe $(M,x_0)$.
La d\'efinition suivante est fondamentale.
\bd
\label{D1}
Une {\em relation ab\'elienne} de la famille $\cal{T}$ est 
un $d$-uplet $\{\phi_1,\ldots,\phi_d\}$,   o\`u $\phi_j$ est une normale ferm\'ee
de~$\cal{F}_j$, dont la somme $\phi_1+\cdots+ \phi_d$ est la forme nulle.
 
Elle est {\em compl\`ete} si aucune des composantes $\phi_j$ n'est la forme nulle.

Le {\em rang} de $\cal T$ est la dimension de l'espace vectoriel
sur $\C$ de ses relations ab\'eliennes.
\ed
La d\'efinition suivante est due \`a Chern et Griffiths \cite{C-Gr3}.
\bd
\label{def-PG}
Une famille de $d$ feuilletages $\cal{F}_j$ de type $(r,n)$ sur
le germe $(M,x_0)$ est {\em un $d$-tissu de type (r,n)}\,  si
ses normales g\'en\'eratrices $\O_j$ v\'erifient la condition 
de position g\'en\'erale (PG) suivante :
$$
\text{(PG)}  \qquad
\;\;\;
1 \leq j_1 < \cdots < j_{\delta} \leq \min(d,n)  \; \Rightarrow \; 
\Omega_{j_1}(x_0) \wedge \cdots \wedge \Omega_{j_\delta}(x_0) \neq 0.
$$
\ed
Dans la premi\`ere partie de Chern-Griffiths~\cite{C-Gr3}, les auteurs d\'emontrent 
la borne suivante pour le rang d'un tissu.  
Nous en rappelons la d\'emonstration dans le Chapitre 2.
\bt
Le rang d'un $d$-tissu de type $(r,n)$ est au plus \'egal,
et cette borne est optimale,
au nombre $\rho_{r,n}(d)$ d\'efini par la formule 
\beq
\label{borne}
\rho_{r,n}(d) = \sum_{h=0}^{+\infty} \, { r-1 +h \choose r-1 } \times \max\,(d-(r+h)(n-1)-1,0).
\eeq
\et
\bd
Un $d$-tissu $\cal{T}$ de type $(r,n)$ est {\em de rang maximal} s'il est de rang $\rho_{r,n}(d)$.
\ed
La borne (\ref{borne}) et donc la d\'efinition d'un tissu de rang 
maximal d\'ependent  bien s\^ur du choix de la condition (PG) dite, mais c'est un peu trompeur,
de position g\'en\'erale.

Sauf dans le cas $r=1$ que nous \'ecartons dans  ce commentaire,
le choix de la condition (PG) n'est pas anodin. 
On pourrait par exemple exiger 
que les normales g\'en\'eratrices $\O_j(x_0)$ des feuilletages
de la famille $\cal{T}$ en $x_0$ 
soient en position g\'en\'erale dans  
l'espace vectoriel $\wedge^rT^\star_{x_0}M$. Ce serait une condition 
plus forte que la condition (PG) si $r\geq 2$.
Pour $d$ assez grand, la borne $\rho_{r,n}(d)$ 
serait remplac\'ee par une borne (beaucoup) 
plus petite et il n'y aurait pas de th\'eor\`eme sur
la classification des tissus de rang maximal. 
La raison du choix de la condition (PG) dans \cite{C-Gr3}
appara\^{\i}tra plus clairement dans 
la section suivante, voir la Remarque \ref{surPG}.

\sk
L'entier suivant interviendra souvent dans cet article :
\beq
\label{q(d)}
q(d) = d - r(n-1) - 2.
\eeq
La seconde partie de Chern-Griffiths~\cite{C-Gr3}  concerne la
g\'eom\'etrie de la famille des normales d'un 
tissu de rang maximal dans le cas $r=2$ et si $q(d)$ est assez grand.
Toutefois les auteurs n'obtiennent un r\'esultat d\'efinitif
que si $n=2$ ou en renfor\c{c}ant l'hypoth\`ese (PG).

\sk
Le Chapitre 2 est consacr\'e \`a ce probl\`eme. Nous 
y \'etablirons la propri\'et\'e attendue de la famille
des normales sous la seule condition $q(d)\geq n-1$,
sans hypoth\`ese sur $(r,n)$ ni 
renforcement de la condition (PG). 

\sk
Nous montrerons en particulier 
qu'un $d$-tissu de rang maximal sur un germe $(M,x_0)$ d\'efinit 
sur ce germe une {\em structure presque grassmannienne} 
si $q(d)\geq n-1$. Il s'agit d'un type de ${\rm G}$-structure
model\'e sur la g\'eom\'etrie de la grassmannienne
des $(n-1)$-plans de $\P^{r+n-1}$.

Pour ne pas alourdir cette introduction, nous en reportons 
la d\'efinition classique au d\'ebut du Chapitre 3.

\ssct{Le th\'eor\`eme d'Abel et sa r\'eciproque}

Soit $\G_{r,n}$ la grassmannienne 
des $(n-1)$-plans de $\P^{r+n-1}$. Comme cette notation, tr\`es commode ici,
n'est pas standard, nous la mettons en \'evidence :
\beq
\label{Grn}
\text{$\G_{r,n}$ est la grassmannienne des $(n-1)$-plans de $\P^{r+n-1}$}. 
\eeq
Pour $x\in \G_{r,n}$, nous notons $H(x)$
cet \'el\'ement vu comme une partie de $\P^{r+n-1}$.
Si $p$ est un point de $\P^{r+n-1}$,
nous notons $\G_{r,n}(p)$ le \og cycle de Schubert \fg\, des 
$(n-1)$-plans qui passent par le point $p$. C'est une sous-vari\'et\'e lisse 
de codimension $r$ de $\G_{r,n}$. Nous appelons ces vari\'et\'es
les {\em sous-vari\'et\'es distingu\'ees} de la grassmannienne.
\bd
Un feuilletage de type $(r,n)$ sur un germe $(\G_{r,n},x_0)$
est grassmannien si ses feuilles sont des ouverts 
de sous-vari\'et\'es distingu\'ees. Une famille de feuilletages 
sur $(\G_{r,n},x_0)$ est {\em une famille grassmannienne} 
si ses membres  sont des feuilletages grassmanniens.
\ed 
Soit $H(x_0)$ un $(n-1)$-plan.
Si l'on note $[\eta_1 : \cdots : \eta_r : \xi_1 : \cdots : \xi_n]$
les coordonn\'ees homog\`enes d'un point, \`a un automorphisme 
de $\P^{r+n-1}$ pr\`es, un $(n-1)$-plan $H(x)$ 
voisin de $H(x_0)$ est d\'efini par un syst\`eme de la forme 
$$
\eta_a - \sum_{\a=1}^n x_{a,\a} \xi_\a = 0 \qquad a=1,\ldots,r.
$$
Les $x_{a,\a}$ param\`etrent localement la grassmannienne. 
Si l'on fixe le  point $p=[\eta_1 : \cdots : \eta_r : \xi_1 : \cdots : \xi_n]$,
le syst\`eme pr\'ec\'edent devient le syst\`eme d'\'equations qui 
d\'efinit la sous-vari\'et\'e distingu\'ee $\G_{r,n}(p)$ de $\G_{r,n}$.

Au voisinage de $x_0$, les sous-vari\'et\'es distingu\'ees sont donc les 
vari\'et\'es int\'egrales des syst\`emes de la forme suivante, o\`u $\xi = [\xi_1 : \cdots : \xi_n]\in \P^{n-1}$ : 
$$
\sum_{\a=1}^n \xi_\a dx_{a,\a}=0, \qquad a=1,\ldots,r.
$$
Une normale en un point $x$ d'un feuilletage grassmannien 
est donc proportionnelle \`a $\wedge_{a=1}^r (\sum_{\a=1}^n \xi_\a dx_{a,\a})$ :
elle appartient au sous-espace {\em de dimension $rn$}
de $\wedge^r T^\star_{x}\G_{r,n}$ engendr\'e
par les $dx_{1,\a_1}\wedge \cdots \wedge dx_{r,\a_r}$.
\bre
\label{surPG}
C'est une cl\'e du choix de la condition (PG) 
dans la D\'efinition \ref{def-PG}. Les normales en un point $x$
d'un tissu grassmannien ne sont jamais en position g\'en\'erale
dans $\wedge^r T^\star_{x}\G_{r,n}$ 
si $r$ est $\geq 2$ et si $d$ est assez grand, tandis que le probl\`eme traditionnel de la th\'eorie 
est de montrer qu'un $d$-tissu de rang maximal est 
isomorphe \`a un tissu grassmannien.
\ere
La notion suivante sera utile dans la suite.
\bd
Une famille finie  de feuilletages $\cal{F}_j$ de type $(r,n)$ 
sur un germe $(M,x_0)$, de normales g\'en\'eratrices $\O_j$,
 est un {\em pr\'e-tissu} si  $\O_j(x_0)\wedge \O_k(x_0)$ est non nul
pour $j\neq k$.
\ed

\bk
Soit $(Z,p)$ un germe de vari\'et\'e lisse de dimension $r$ 
en $p\in \P^{r+n-1}$, transverse au $(n-1)$-plan $H(x_0)$.
Un $(n-1)$-plan voisin $H(x)$ coupe $(Z,p)$ en un point 
$\k(x)$ voisin de $p$. Le morphisme d'incidence 
$$
\k: (\G_{r,n},x_0) \rightarrow (Z,p), \qquad x\mapsto \k(x),
$$
est une submersion. Il d\'efinit un feuilletage  
$\cal{F}$ sur $(\G_{r,n},x_0)$, dont les feuilles sont les fibres de $\k$.
C'est un feuilletage grassmannien et l'on v\'erifie facilement
que tout feuilletage grassmannien s'obtient ainsi.

Le feuilletage $\cal{F}$ est le {\em feuilletage d'incidence} 
d\'efini par le germe $(Z,p)$.

\sk
Le lemme suivant est tautologique mais il explique la d\'efinition 
initiale par Blaschke d'une relation ab\'elienne.
\ble
\label{tauto}
L'application $\phi\mapsto \k^\star \phi$ induit un isomorphisme 
de l'espace des $r$-formes lisses sur le germe $(Z,p)$ sur l'espace des 
normales ferm\'ees du feuilletage d'incidence $\cal{F}$ sur $(\G_{r,n},x_0)$ 
qu'il d\'efinit.
\ele
\bpf
Il est clair par d\'efinition que $\k^\star \phi$ 
est une normale de $\cal{F}$ si $\phi$
est une $r$-forme lisse sur $(Z,p)$.
C'est une normale ferm\'ee puisque $\phi$ 
est de degr\'e maximal $r$ donc ferm\'ee.
R\'eciproquemant, si $\phi$ ne s'annule pas,
$\k^\star \phi$ est une normale g\'en\'eratrice ferm\'ee 
de $\cal{F}$ et une normale g\'en\'erale $f\k^\star\phi$ est ferm\'ee 
si et seulement si $df\wedge \k^\star \phi=0$,
si et seulement si $f$ est de la forme $u\circ \k$.
Alors $f\k^\star \phi= \k^\star (u\phi)$.
\epf
Pour obtenir un pr\'e-tissu (respectivement un tissu)
d'incidence, ou grassmannien 
puisque c'est la m\^eme chose, d'ordre $d$ sur le germe 
$(\G_{r,n},x_0)$, il faut et il suffit de se donner
$d$ germes $(Z_j,p_j)$ de vari\'et\'es lisses
de dimension~$r$, transverses \`a $H(x_0)$ 
en $d$ points $p_j$ deux-\`a-deux distincts
(respectivement en position g\'en\'erale dans $H(x_0)$).
Nous avons alors $d$ morphismes d'incidence 
$$
\k_j: (\G_{r,n},x_0) \rightarrow (Z_j,p_j).
$$
Compte tenu du lemme pr\'ec\'edent,
l'application 
$$
(\phi_1,\ldots,\phi_d)\mapsto (\k_1^\star \phi_1, \ldots, \k_d^\star \phi_d),
$$
o\`u les $\phi_j$ sont des $r$-formes lisses sur les germes $(Z_j,p_j)$,
induit un isomorphisme de l'espace des $d$-uplets 
de {\em trace nulle} sur l'espace des relations ab\'eliennes du pr\'e-tissu.
Ici la trace est la somme  
$$
\k_1^\star \phi_1 +  \cdots + \k_d^\star \phi_d.
$$

\sk
L'exemple o\`u les $(Z_j,p_j)$ sont des germes d'une m\^eme 
vari\'et\'e projective est fondamental.
\bd
Soit $Z$ une sous-vari\'et\'e projective de $\P^{r+n-1}$,
r\'eduite de dimension pure $r$ et de degr\'e $d$. Une {\em $r$-forme ab\'elienne}
sur $Z$ est une $r$-forme d\'efinie et lisse
sur un ouvert dense de $Z$ dont la trace, d\'efinie 
au voisinage du point g\'en\'eral de $\G_{r,n}$, est nulle.
Le {\em genre ab\'elien} de $Z$ est la dimension de l'espace 
de ses $r$-formes ab\'eliennes.
\ed
Plus pr\'ecis\'ement, $Z$ coupe un $(n-1)$-plan g\'en\'eral $H(x)$
transversalement 
en $d$ points deux-\`a-deux distincts $p_j$ et la forme 
$\phi$ consid\'er\'ee est d\'efinie au voisinage des $p_j$
dans $Z$. Au voisinage de $x$, nous avons comme ci-dessus $d$
morphismes d'incidence $\k_j$. La trace de $\phi$
est d\'efinie localement comme \'etant la somme $\k_1^\star \phi +  \cdots + \k_d^\star \phi$.

\bk
Tout ceci a bien s\^ur \`a voir avec le th\'eor\`eme
d'addition d'Abel et ses g\'en\'eralisations.
Dans \cite{Gr} Griffiths montre que les formes ab\'eliennes,
qu'il nomme autrement, sont rationnelles sur $Z$ et lisses sur la partie r\'eguli\`ere de $Z$.
Il les caract\'erise
aussi lorsque $Z$ est une hypersurface. 
Il d\'ecoule de ses r\'esultats 
que le genre ab\'elien de $Z$ est \'egal 
\`a son genre g\'eom\'etrique si $Z$ est lisse
et \`a son genre arithm\'etique si $Z$ est une hypersurface.

\sk
Il a \'et\'e remarqu\'e plus r\'ecemment\footnote{ Voir
Henkin-Passare \cite{H-P}. Signalons que  
les travaux que nous citons dans cette section
concernent des formes de degr\'e quelconque, pas seulement de degr\'e 
maximal.} que 
les $r$-formes ab\'eliennes sur $Z$ sont les sections 
globales du faisceau des $r$-formes de Barlet, voir Barlet~\cite{Ba}
pour les d\'efinitions et la d\'emonstration, en particulier 
la deuxi\`eme partie de cet article dans laquelle 
l'auteur \'etablit ce qu'il appelle la propri\'et\'e de la trace 
universelle.
Dans le cas qui nous occupe, le faisceau de Barlet 
est un faisceau de $r$-formes rationnelles, 
intrins\`equement attach\'e \`a la vari\'et\'e {\em abstraite} $Z$,
et isomorphe au faisceau $\o^r_Z$, 
le faisceau dualisant de Grothendieck de $Z$.
Ainsi : 
$$
\{ \text{\small $r$-formes ab\'eliennes sur $Z$} \}
=
\{ \text{\small $r$-formes de Barlet sur $Z$} \}
\simeq 
H^0(Z,\o^r_Z).
$$

Il en r\'esulte, et c'est remarquable, que les formes ab\'eliennes
sur $Z$ et le genre ab\'elien de $Z$ ne d\'ependent 
pas de la repr\'esentation de $Z$ comme sous-vari\'et\'e d'un espace projectif.
Comme ces consid\'erations, qui sont importantes d'un point de vue th\'eorique, 
ne joueront pas de r\^ole dans cet article, nous ne donnons pas plus de d\'etail.
\bd
\label{Conv1}
{\em Par abus de langage nous dirons qu'une sous-vari\'et\'e $Z$
de $\P^{r+n-1}$, r\'eduite de dimension pure $r$ et de degr\'e $d$,
d\'efinit {\em un pr\'e-tissu 
alg\'ebrique d'incidence ou grassmannien} $\cal{T}_Z$ d'ordre $d$ sur $\G_{r,n}$, 
dont le pr\'e-tissu qu'elle d\'efinit en un point g\'en\'eral
est un germe.}
\ed
La raison de cet abus de langage est qu'\`a notre connaissance,
la notion globale de \og (pr\'e-)tissu alg\'ebrique singulier \fg\, 
n'a pas \'et\'e d\'efinie en toute g\'en\'eralit\'e.
Il para\^{\i}t toutefois \'evident que l'objet $\cal{T}_Z$,
d\'efini ci-dessus comme une collection de germes de tissus,
m\'eritera d'\^etre appel\'e, le jour venu, un pr\'e-tissu alg\'ebrique, 
en l'occurence d'un type tr\`es particulier.

\bk
La r\'eciproque suivante du th\'eor\`eme d'Abel ou {\em th\'eor\`eme d'Abel inverse}
est essentielle. L'id\'ee  de la d\'emonstration remonte \`a Darboux~\cite{Da}.
Le r\'esultat g\'en\'eral est d\^u \`a Griffiths~\cite{Gr}.
\bt
\label{abel-inv}
Soit $\cal{T}$ un pr\'e-tissu grassmannien de type $(r,n)$ et d'ordre $d$ 
sur $(\G_{r,n},x_0)$, d\'efini par $d$ germes lisses
$(Z_j,p_j)$ de dimension~$r$, transverses au $(n-1)$-plan $H(x_0)$
en des points diff\'erents.

Si $\cal{T}$ admet une relation 
ab\'elienne {\em compl\`ete}, $\cal{T}$ est un germe du pr\'e-tissu alg\'ebrique 
grassmannien d\'efini par une vari\'et\'e $Z$ de degr\'e $d$ qui 
contient les $(Z_j,p_j)$. Le rang de $\cal{T}$ est 
\'egal au genre ab\'elien de $Z$ et ses relations ab\'eliennes sont 
induites par les $r$-formes ab\'eliennes sur $Z$.
\et

\ssct{M\'ethode canonique et vari\'et\'e de Blaschke} 

Il s'agit de d\'ecrire la {\em m\'ethode canonique},  d\'evelopp\'ee par Blaschke, Bol et Howe
\`a partir d'une id\'ee de Poincar\'e \cite{Po},
puis  de d\'efinir la {\em vari\'et\'e de Blaschke} d'un tissu de rang 
maximal, qui intervient dans l'\'enonc\'e de nos r\'esultats.

\sk
Soit $\cal{T}=\{\cal{F}_1,\ldots,\cal{F}_d\}$ un $d$-tissu 
de type $(r,n)$ sur un germe lisse $(M,x_0)$. Notons $N+1$
le rang de $\cal{T}$. La m\'ethode consiste \`a se donner une base 
$$
\phi^{(\l)} = (\phi_1^{(\l)}, \ldots ,\phi^{(\l)}_d), \qquad \l=1,\ldots,N+1,
$$
de l'espace des relations ab\'eliennes du tissu $\cal{T}$ et \`a d\'efinir 
les {\em applications de Poincar\'e} $\k_j$, une pour chaque feuilletage $\cal{F}_j$, par :
$$
\k_j: (M,x_0)\dasharrow \P^N, \qquad \k_j(x) = [\phi^{(1)}_j(x): \cdots :\phi^{(N+1)}_j(x)].
$$
Pour que ces applications soient d\'efinies comme applications m\'ero\-morphes,
nous  supposons  que le tissu $\cal{T}$ admet au moins une relation 
ab\'elienne compl\`ete.

\bk
Avant de continuer et pour motiver la construction pr\'ec\'edente, il 
est certainement utile de montrer sa relation avec la construction 
de la {\em vari\'et\'e canonique} d'une vari\'et\'e projective
en g\'eom\'etrie alg\'ebrique. Pour la d\'efinir, nous utilisons 
les $r$-formes ab\'eliennes plut\^ot que les $r$-formes 
lisses sur $Z$. On retrouve la d\'efinition la plus classique si $Z$
est lisse.

Cette relation appara\^{\i}t si l'on choisit comme tissu $\cal{T}$
un germe d'un tissu alg\'ebrique grassmannien
$\cal{T}_Z$, d\'efini par une sous-vari\'et\'e alg\'ebrique $Z$ 
de $\P^{r+n-1}$, de dimension $r$ et de degr\'e $d$, de genre ab\'elien 
$N+1$. Nous 
supposons qu'il existe une $r$-forme ab\'elienne 
sur $Z$ qui ne s'annule sur aucune des composantes irr\'eductibles de $Z$.
Choisissons une base $(\phi^{(1)},\ldots,\phi^{(N+1)})$ 
de l'espace des $r$-formes ab\'eliennes sur $Z$ et consid\'erons l'application
$$
\k_Z : \; Z \dasharrow \P^N, \qquad \k_Z(p) = [\phi^{(1)}(p) : \ldots : \phi^{(N+1)}(p)].
$$
L'image de cette application est {\em la vari\'et\'e 
canonique de $Z$}. Elle est ind\'ependante du choix de la base, \`a un
automorphisme de $\P^N$ pr\`es.

\sk
Soit alors $H(x_0)$ un $(n-1)$-plan coupant $Z$ transversalement en $d$ points 
lisses $p_j$ de $Z$ et $\chi_j: (\G_{r,n},x_0)\rightarrow (Z,p_j)$
les morphismes d'incidence.
Les $d$-uplets $(\chi_1^\star \phi^{(\l)},\ldots, \chi_d^\star\phi^{(\l)})$
forment une base de l'espace des relations ab\'eliennes du germe 
en $x_0$ du tissu $\cal{T}_Z$ et pour ce choix de base, les applications de 
Poincar\'e $\k_j$ sont donn\'ees par 
$$
x\in (M,x_0), \qquad \k_j(x) = [\phi^{(1)}(\chi_j(x)) : \cdots : \phi^{(N+1)}(\chi_j(x))].
$$

\sk
{\em Ceci montre que les images des applications de Poincar\'e associ\'ees 
\`a un germe du tissu alg\'ebrique $\cal{T}_Z$
sont des germes de la vari\'et\'e canonique  de $Z$, 
d\'efinie \`a partir d'une base de $r$-formes ab\'eliennes sur $Z$.}

\bk
Nous revenons \`a la m\'ethode canonique appliqu\'ee \`a un tissu.
Rappelons que l'entier $q(d)$ est d\'efini par (\ref{q(d)}).
La m\'ethode canonique ne semble pas pouvoir donner de r\'esultat 
int\'eressant si $q(d)<n-1$. En revanche, si $q(d)\geq n-1$, elle nous en 
donnera un sous une hypoth\`ese plus faible que celle pour un 
tissu d'\^etre de rang maximal. La notion suivante 
a \'et\'e introduite dans Tr\'epreau \cite{Tr} sous le 
nom plus \'evocateur mais trop long de \og tissu de rang maximal en 
valuation $\leq 1$ \fg.
\bd
\label{semi}
Un $d$-tissu de type $(r,n)$ sur un germe $(M,x_0)$
est {\em  semi-extr\'emal} si $q(d)\geq n-1$ et s'il admet 
$$
[d-r(n-1)-1] + r[d - (r+1)(n-1)-1]
$$
relations ab\'eliennes dont les $1$-jets en $x_0$ sont
lin\'eairement ind\'ependants.
\ed
Ce nombre est en fait le plus grand possible tel qu'on puisse avoir la 
propri\'et\'e pr\'ec\'edente si $q(d)\geq n-1$. De plus 
un $d$-tissu de rang maximal est semi-extr\'emal si $q(d)\geq n-1$,
voir le Chapitre 3.

\sk
Les propri\'et\'es g\'eom\'etriques du syst\`eme des normales d'un tissu 
semi-extr\'emal, \'etablies dans le Chapitre 2, nous permettront dans le Chapitre~3
d'obtenir les propri\'et\'es suivantes,  
d\'ej\`a connues si $r=1$, pour les tissus de codimension $r\geq 2$.

D'abord les applications $\k_j$ sont de rang $r$. L'image de $\k_j$
est donc un germe lisse $(\tilde{Z}_j,\k_j(x_0))$ de dimension $r$
et les fibres de la submersion induite 
$$
\k_j: (M,x_0)\rightarrow (\tilde{Z}_j,\k_j(x_0))
$$
sont les feuilles du feuilletage $\cal{F}_j$.
Ensuite, si l'on choisit $n$ feuilletages parmi les $\cal{F}_j$, 
par exemple $\cal{F}_1,\ldots, \cal{F}_n$, l'application
$$
(M,x_0) \rightarrow \prod_{j=1}^n (\tilde{Z}_j,\k_j(x_0)), \qquad x\mapsto (\k_1(x),\ldots,\k_n(x)),
$$
est un isomorphisme. 
Enfin, pour $x\in (M,x_0)$, les points ${\k}_j(x)$, $j$ variant de $1$ \`a $d$,
 sont deux-\`a-deux distincts et situ\'es 
sur une et une seule m\^eme courbe rationnelle normale $\s(x)$ de degr\'e $q(d)$. 
Il suffit ici de savoir qu'on appelle ainsi une courbe rationnelle irr\'eductible
et lisse de degr\'e $q(d)$, qui engendre un espace projectif de dimension $q(d)$. 
\bd
\label{var-X}
{\em La vari\'et\'e de Blaschke} $X_{\cal{T}}$ du tissu semi-extr\'emal~$\cal{T}$
est l'intersection de toutes les sous-vari\'et\'es projectives de $\P^N$
qui contiennent les courbes rationnelles normales $\s(x)$, $\,x\in (M,x_0)$.
\ed
Il est clair que {\em la famille} 
des germes $(\tilde{Z}_j,\k_j(x_0))$ est ind\'ependante
du choix de la base de relations ab\'eliennes \`a un automorphisme de $\P^N$ pr\`es.
Il en va donc de m\^eme de la vari\'et\'e $X_{\cal{T}}$.
{\em Ce sont des invariants projectifs  
de la classe d'isomorphie du tissu $\cal{T}$.}

\sk
La m\'ethode canonique permet ainsi d'associer des objets de nature g\'eom\'etrique
\`a un tissu semi-extr\'emal et offre une 
strat\'egie pour montrer qu'il est isomorphe \`a un germe de tissu 
grassmannien : montrer que les 
$(\tilde{Z}_j,\k_j(x_0))$ sont des germes d'une 
vari\'et\'e projective $\tilde{Z}$ et que celle-ci est la vari\'et\'e 
canonique d'une vari\'et\'e projective $Z$.

\ssct{R\'esultats d'alg\'ebricit\'e}

Nous \'enon\c{c}ons dans cette section les principaux r\'esultats 
de cet article. Sauf pour le premier, leurs d\'emonstrations dans 
le Chapitre 4 repose aussi sur le travail 
pr\'eliminaire accompli dans Pirio-Tr\'epreau \cite{P-Tr}.
Le premier r\'esultat concerne la vari\'et\'e de Blaschke d'un tissu 
semi-extr\'emal. 
\bt
\label{Th-X0}
La vari\'et\'e de Blaschke $X_{\cal{T}}$ d'un $d$-tissu semi-extr\'emal $\cal{T}$, 
de type $(r,n)$ avec $r\geq 2$,  est une vari\'et\'e projective 
irr\'eductible de dimension
$r+1$. De plus, pour  $(a_1,\ldots,a_n)\in X_{\cal{T}}^n$ g\'en\'erique, il existe une courbe
rationnelle normale de degr\'e $q(d)$, contenue dans $X_\cal{T}$ et passant par
les points $a_1,\ldots,a_n$.
\et
Le point cl\'e est que la vari\'et\'e $X_{\cal{T}}$ est de dimension
$r+1$. Le reste en sera une cons\'equence facile,
compte tenu des propri\'et\'es des applications de Poincar\'e
que nous avons mentionn\'ees plus haut.

Le r\'esultat est d\^u \`a H\'enaut \cite{He} pour $n=2$ et \`a Tr\'epreau \cite{Tr}
pour $r=1$ avec $n\geq 3$. La d\'emonstration pour $r=1$ et $n\geq 3$ est plus difficile 
et ne sera pas reproduite ici. Le calcul pour $r\geq 2$ reste encore un peu compliqu\'e 
et donnera lieu aux quelques pages techniques de cet article.
\bre
Rappelons \`a quoi correspond, au moins dans certains cas, la vari\'et\'e de Blaschke 
du tissu $\cal{T}_Z$ d\'efini par une vari\'et\'e projective $Z\subset \P^{r+n-1}$,
de dimension pure $r$, de degr\'e $d$ avec $q(d)\geq n-1$, 
et qui rencontre un $(n-1)$-plan g\'en\'erique en $d$ points 
en position g\'en\'erale dans ce plan, de genre ab\'elien 
maximal $N+1=\rho_{r,n}(d)$. 

\sk
{\em Si $Z$ est lisse} (voir Harris \cite{Ha}), la vari\'et\'e  
$Z$ est contenue dans une sous-vari\'et\'e $X_0$ de $\P^{r+n-1}$, de dimension
$r+1$ et de {\em degr\'e minimal} $\,n-1$. De plus, 
l'application canonique $\k_Z: Z \dasharrow \P^N$
admet un prolongement  naturel \`a $X_0$. L'image de $X_0$
par ce prolongement n'est autre que la vari\'et\'e de Blaschke du tissu $\cal{T}_Z$.

Il en va probablement de m\^eme sous la seule 
hypoth\`ese que le tissu $\cal{T}_Z$ est semi-extr\'emal
mais nous ne le d\'emontrons pas ici.
\ere
Le Th\'eor\`eme \ref{Th-X0} est une \'etape importante 
dans l'application de la m\'ethode canonique. Comme on verra, il permet
de construire un mod\`ele g\'eom\'etrique d'un tissu semi-extr\'emal,
mais il reste \`a montrer que ce mod\`ele est alg\'ebrique.

\sk
Le moment est propice pour situer le pr\'esent article par rapport \`a 
Tr\'epreau \cite{Tr}.
Le th\'eor\`eme pr\'ec\'edent est vrai aussi si $r=1$ avec $n\geq 3$
et permet alors de d\'emontrer qu'un $d$-tissu $\cal{T}$ semi-extr\'emal de type
$(1,n)$ avec $n\geq 3$ est isomorphe \`a un germe d'un  tissu 
alg\'ebrique grassmannien. Dans ce cas, la vari\'et\'e de Blaschke $X_\cal{T}$
{\em est une surface} et les courbes rationnelles normales 
de degr\'e $q(d)$ contenues dans $X_{\cal{T}}$ sont des diviseurs.
Cette particularit\'e permet d'appliquer, comme le faisait Bol~\cite{BB}, un r\'esultat classique 
d'Enriques sur la lin\'earit\'e de certains 
syst\`emes alg\'ebriques de diviseurs et d'en d\'eduire directement 
que le tissu $\cal{T}$ est isomorphe \`a un germe de tissu 
alg\'ebrique grassmannien.

\sk
Ceci n'est plus possible si $r\geq 2$. Alors $X_\cal{T}$ est de dimension $\geq 3$,
les courbes ne sont pas des diviseurs et nous ne savons pas conclure
\`a partir des seules propri\'et\'es de $X_\cal{T}$ donn\'ees par le Th\'eor\`eme \ref{Th-X0}, 
m\^eme dans la cas le plus simple $n=2$,
sans faire l'hypoth\`ese que le tissu est
de rang maximal. De plus, sous cette hypoth\`ese, la solution s'est av\'er\'ee  
beaucoup plus difficile \`a obtenir que dans \cite{Tr}. 

\bk
Dans la perspective de surmonter ces difficult\'es, 
nous avons introduit et \'etudi\'e dans Pirio-Tr\'epreau \cite{P-Tr} les 
classes de vari\'et\'es suivantes.
\bd
\label{r-n-q}
Soit $q\geq n-1$ et posons $d=q+r(n-1)+2$.
Une sous-vari\'et\'e  projective $X$ de~$\P^N$,
irr\'eductible de dimension $r+1$ et non contenue dans un hyperplan 
de~$\P^N$, appartient \`a la classe 
$\cal{X}_{r+1,n}(q)$ si  les deux conditions suivantes 
sont satisfaites :
\be
\item[(a)] pour $(a_1,\ldots,a_n)\in X^n$ g\'en\'erique, il existe une courbe
rationnelle normale de degr\'e $q$, contenue dans $X$ et passant par
$a_1,\ldots,a_n$ ;

\item[(b)]  la dimension $N$ v\'erifie 
$N+1=\rho_{r,n}(d)$\footnote{En fait $N+1$ est donn\'e sous une 
autre forme dans \cite{P-Tr}. Nous v\'erifierons 
dans le Chapitre 4 qu'on a bien $N+1=\rho_{r,n}(d)$.}.
\ee
\ed
La condition (a) seule implique 
qu'on a $N+1\leq \rho_{r,n}(d)$, voir \cite{P-Tr}~Th\'eor\`eme 1.2 :
la propri\'et\'e (b) signifie que $X$ engendre un espace de dimension
maximale, compte tenu de la propri\'et\'e (a).

\bk
La vari\'et\'e de Blaschke d'un $d$-tissu {\em de rang maximal} de type 
$(r,n)$ appartient donc \`a la classe $\cal{X}_{r+1,n}(q(d))$
si $r\geq 2$ et $q(d)\geq n-1$. 

\bk
Nous en venons maintenant aux principaux r\'esultats de cet article.
\bd
Deux $d$-tissus sur des 
germes $(M,x_0)$ et $(M',x'_0)$
sont {\em isomorphes} s'il existe un 
isomorphisme de $(M,x_0)$ sur $(M',x'_0)$ qui les \'echangent.
\ed
Le probl\`eme traditionnel est de montrer qu'un $d$-tissu $\cal{T}$
de rang maximal et de type $(r,n)$, nous supposons ici $r\geq 2$, 
est isomorphe \`a un germe d'un tissu alg\'ebrique grassmannien
si $d$ est assez grand.

\sk
Nous avons le r\'esultat suivant qui r\'esout en particulier le probl\`eme 
pr\'ec\'edent. Nous le donnons en 
premier parce que toutes les notions qui interviennent dans l'\'enonc\'e 
ont d\'ej\`a \'et\'e introduites. Rappelons la notation $q(d) = d - r(n-1) - 2$.

\bt
\label{Th1A}
Soit $\cal{T}$ un $d$-tissu de type $(r,n)$ et de rang maximal, avec $r\geq 2$ et $q(d)\geq n-1$.
Si $n=2$ ou si $q(d)\neq 2n-3$, il est isomorphe \`a un germe d'un tissu 
alg\'ebrique grassmannien.
\et
Avant de le commenter, disons qu'il d\'ecoule de l'un des r\'esultats essentiels 
de  \cite{P-Tr}, la classification compl\`ete des vari\'et\'es 
de la classe $\cal{X}_{r+1,n}(q)$, pr\'ecis\'ement  sous l'une des deux conditions 
$n=2$ ou $q\neq 2n-3$. Nous reviendrons plus bas sur le cas 
$q=2n-3$ avec $n\geq 3$, exclu de l'\'enonc\'e.

\sk
Concernant l'hypoth\`ese $q(d)\geq n-1$, il semble qu'elle soit 
n\'ecessaire pour que la m\'ethode canonique puisse 
\^etre appliqu\'ee avec succ\`es. 
Comme le rang de $\cal{T}$ est nul si $q(d)<0$, il faut
bien s\^ur supposer $q(d)\geq 0$. La question int\'eressante
de l'existence ou non de contre-exemples lorsque 
$0\leq q(d) < n-1$ est ouverte en g\'en\'eral.
Les seuls contre-exemples connus concernent des tissus de 
codimension~$2$ avec $q(d)=0$,  voir en particulier Goldberg \cite{Go1}
pour $(r,n)=(2,2)$.

\sk
Notre r\'esultat est compl\`etement nouveau sauf si $q(d)=n-1$ ou si $n=2$.

\sk
On peut dire que l'\'enonc\'e ci-dessus pour  $q(d)=n-1$ est folklorique.
La d\'emonstration de Howe (1932) lorsque $(r,n)=(1,2)$ 
se g\'en\'eralise sans encombre. Nous y reviendrons dans 
le Chapitre 3.

\sk
Le cas $n=2$ a \'et\'e \'etudi\'e par H\'enaut \cite{He}, qui suit la m\'ethode canonique,
mais sa d\'emonstration est d\'efaillante sur deux points.
D'abord, concernant la g\'eom\'etrie des normales d'un tissu 
de type $(r,2)$,  l'auteur se r\'ef\`ere \`a un \'enonc\'e 
de Little \cite{Li}, correct dans ce cas mais dont la d\'emonstration 
dans \cite{Li} est au moins douteuse, voir le Chapitre 2. 
D'autre part, H\'enaut \'enonce bien que, sous les hypoth\`eses 
du Th\'eor\`eme~\ref{Th1A} avec $n=2$, 
la vari\'et\'e de Blaschke du tissu appartient \`a ce que nous appelons 
la classe $\cal{X}_{r+1,2}(q(d))$ mais son argument pour en d\'eduire que 
c'est une vari\'et\'e de Veronese  d'ordre $q(d)$, ce qui est n\'ecessaire pour conclure,  
est incorrect, voir le Chapitre 4.

\bk
Rappelons que l'\'enonc\'e pr\'ec\'edent ne concerne que les 
tissus {\em de rang maximal} avec $q(d)\geq n-1$. 
Si de plus $n=2$ ou si $q(d)\neq 2n-3$, il caract\'erise 
donc les $d$-tissus de type $(r,n)$ qui sont isomorphes \`a un 
germe d'un tissu $\cal{T}_Z$, o\`u $Z\subset \P^{r+n-1}$ 
est de dimension pure $r$, de degr\'e $d$, rencontre 
un $(n-1)$-plan g\'en\'erique en $d$ points en position g\'en\'erale dans ce plan
et {\em est de genre ab\'elien maximal} compte tenu de ces conditions.

\sk
Le cas $n=2$ a ceci de particulier que le genre ab\'elien d'une hypersurface d'un espace 
projectif est toujours maximal puisqu'il est \'egal \`a son 
genre arithm\'etique (Griffiths, voir la Section 1.2) et ne d\'epend donc que de son degr\'e.
On a donc :
\bco
\label{Th1D}
Si $r\geq 2$ et $d\geq r+3$, un $d$-tissu de type $(r,2)$ est de rang maximal si et seulement s'il est 
isomorphe \`a un germe d'un tissu alg\'ebrique d\'efini
par une hypersurface de degr\'e $d$ de $\P^{r+1}$.
\eco

\sk
Le prochain \'enonc\'e donnera  une repr\'esentation alg\'ebrique 
canonique de tout $d$-tissu de rang maximal de type $(r,n)$, toujours  avec $r\geq 2$ et $q(d)\geq n-1$, 
mais sans autre 
hypoth\`ese. Sa compr\'ehension exige que nous fassions quelques rappels tir\'es de \cite{P-Tr}.
Des r\'ef\'erences pr\'ecises seront donn\'ees 
dans le Chapitre 4.

\sk
Soit $X$ une vari\'et\'e de la classe $\cal{X}_{r+1,n}(q)$.
On montre qu'il existe une unique vari\'et\'e projective 
de $1$-cycles effectifs de $X$, nous la notons $\Si_X$,
dont l'\'el\'ement g\'en\'erique est une courbe rationnelle 
de degr\'e $q$ et telle que 
pour tout $(a_1,\ldots,a_n)\in X^n$, il existe  un \'el\'ement 
de $\Si_X$ dont le support dans $X$ contient les points $a_j$.

Pour plus de clart\'e, si $x\in \Si_X$, nous notons $\s(x)$
l'\'el\'ement $x$ vu comme une courbe (\'eventuellement non r\'eduite)
de $X$.

\sk
La vari\'et\'e $\Si_X$ est irr\'eductible de dimension $rn$.
On d\'efinit aussi un ouvert Zariski-dense $\Si_{X,{\rm adm}}$
de $\Si_X$, nous rappellerons comment dans le Chapitre 4.
Un \'el\'ement $x$ de $\Si_{X,{\rm adm}}$ est appel\'e 
un {\em \'el\'ement admissible} de $\Si_X$ et l'on dit que $\s(x)$
 est une {\em courbe admissible} de $X$.

On a les propri\'et\'es suivantes. 
Un \'el\'ement admissible $x$
de $\Si_X$ est un point lisse de $\Si_X$ et la courbe 
$\s(x)$ est une courbe rationnelle normale de degr\'e $q$, contenue 
dans la partie lisse de $X$. La r\'eunion de ces courbes 
est un ouvert Zariski-dense $X_{\rm adm}$ de la partie lisse de $X$.

D'autre part, si $\s(x)$ est une courbe admissible et $p$ un point 
de $\s(x)$, l'ensemble alg\'ebrique 
$$
\Si_X(p) = \{x'\in \Si_X, \;\; p\in \s(x')\}
$$
est, au voisinage de $x$, lisse de codimension $r$ dans $\Si_X$.
Enfin, si $p_1,\ldots,p_n$ sont des points deux-\`a-deux distincts 
de $\s(x)$, les espaces tangents $T_x\Si_X(p_j)$ sont en position g\'en\'erale 
dans $T_x\Si_X$.

\bk
Nous avons vu dans la Section 1.2 que la paire $(\P^{r+n-1},\G_{r,n})$
permet d'associer \`a une sous-vari\'et\'e alg\'ebrique $Z\subset \P^{r+n-1}$,
r\'eduite de dimension pure $r$ et de degr\'e $d$, un pr\'e-tissu 
alg\'ebrique d'incidence $\cal{T}_Z$ sur $\G_{r,n}$, 
dont le rang est \'egal au genre ab\'elien de $Z$.

Une construction analogue est possible pour la paire $(X,\Si_X)$
et une hypersurface alg\'ebrique $Z$ de $X$, \`a partir des 
propri\'et\'es d'incidence de $Z$ avec un \'el\'ement g\'en\'eral de 
$\Si_X$, comme nous allons voir. Nous faisons 
l'hypoth\`ese suivante.

\sk
{\em L'hypersurface $Z$ est r\'eduite et toutes ses composantes irr\'eductibles 
rencontrent la partie admissible $X_{\rm adm}$ de $X$.}

\sk
La premi\`ere condition nous est famili\`ere.
La seconde condition, qui n'en est une que si $X\bck X_{\rm adm}$ est de codimension un 
dans $X$, n'est pas vraiment n\'ecessaire. Simplement 
une composante de $Z$ contenue dans $X\bck X_{\rm adm}$
ne rencontrerait  aucune courbe admissible de $X$ et ne jouerait 
aucun r\^ole dans la construction qui suit.
Le plus simple est d'\'ecarter ces \'eventuelles composantes.

\sk
Comme une courbe admissible $\s(x)$ est contenue dans $X_{\rm reg}$, 
son nombre d'intersection $\s(x)\cdot Z$ avec $Z$
est bien d\'efini par la th\'eorie \'el\'ementaire 
de l'intersection. Il est ind\'ependant de la courbe admissible 
choisie et c'est un entier $d$ strictement positif puisque $Z$
rencontre $X_{\rm adm}$ par hypoth\`ese.

\sk
Nous v\'erifierons que si la courbe admissible
$\s(x_0)$ est assez g\'en\'erale, elle rencontre $Z$ transversalement 
en $d$ points deux-\`a-deux distincts de $Z_{\rm reg}$.
Nous avons donc, au voisinage de $x_0$, $d$ morphismes d'incidence 
$$
\k_j: (\Si_X,x_0)\rightarrow (Z,p_j),
$$
o\`u $p_j=\k_j(x_0)$ et $\k_j(x)$ est le point d'intersection de $\s(x)$ avec 
$(Z,p_j)$.
Nous v\'erifierons aussi que chaque morphisme $\k_j$
est une submersion. Il d\'efinit donc un feuilletage $\cal{F}_j$
de type $(r,n)$
sur $(\Si_X,x_0)$ dont les feuilles sont les fibres de $\k_j$,
des ouverts au sens usuel de sous-vari\'et\'es de 
la forme $\Si_X(p)$.
Compte tenu de la propri\'et\'e des espaces tangents $T_{x_0}\Si_X(p_j)$
que nous avons mentionn\'ee plus haut,
les feuilletages $\cal{F}_j$ forment un tissu.

Nous nous autorisons le m\^eme abus de langage que dans le cas 
de la paire $(\P^{r+n-1},\G_{r,n})$.
\bd
\label{Conv2}
{\em Nous dirons qu'une hypersurface alg\'ebrique r\'eduite $Z$,
dont toutes les composantes irr\'eductibles rencontrent $X_{\rm adm}$,
d\'efinit sur $\Si_X$ {\em un tissu 
alg\'ebrique d'incidence} $\cal{T}_Z$ pour la paire $(X,\Si_X)$,
d'ordre le nombre d'intersection $d$ de $Z$ avec une courbe admissible 
de $X$, dont le tissu qu'elle d\'efinit au point g\'en\'eral de $\Si_X$
est un germe.}
\ed
Si $\phi_j$ est un germe de $r$-forme lisse sur $(Z,p_j)$,
sa trace $\sum_{j=1}^d \k_j^\star \phi_j$
est un germe de forme sur $(\Si_X,x_0)$ et, dans le cas 
o\`u les $\phi_j$ sont des germes d'une m\^eme 
$r$-forme lisse $\phi$, d\'efinie sur un ouvert Zariski-dense de $Z$,
pour la m\^eme raison que dans la situation analogue 
d\'ecrite dans la Section 1.2, les traces qu'on obtient ainsi sont
des germes d'une forme d\'efinie sur un ouvert dense de $\Si_X$.
Il est clair aussi que l'analogue du Lemme \ref{tauto} reste vrai dans 
ce nouveau cadre. 
\bd
\label{sigma-rel}
Une {\em $r$-forme $\Si_X$-ab\'elienne} sur $Z$ est une $r$-forme lisse sur un ouvert Zariski-dense 
de $Z$, de trace nulle.
\ed
Nous montrerons que le rang du tissu $\cal{T}_Z$ au point g\'en\'eral de $\Si_X$
est \'egal \`a la dimension de l'espace des $r$-formes 
$\Si_X$-ab\'eliennes sur $Z$ et que celles-ci sont rationnelles. 
Nous ne savons pas si elles co\"{\i}ncident toujours 
avec les formes ab\'eliennes, c'est-\`a-dire avec les $r$-formes de Barlet sur $Z$.

\sk
Nous avons le r\'esultat suivant.
\bt
\label{Th1B}
Soit $\cal{T}$ un $d$-tissu de type $(r,n)$ et de rang maximal,  avec $r\geq 2$ et $q(d)\geq n-1$.
Soit $X_{\cal{T}}$ sa vari\'et\'e de Blaschke. 

Il existe une hypersurface projective $Z$ de $X_\cal{T}$,
d\'efinissant un tissu alg\'ebrique d'incidence 
$\cal{T}_Z$, tel que le tissu $\cal{T}$ est isomorphe \`a un germe 
du tissu $\cal{T}_Z$.
\et
Il n'est pas exclu ou plut\^ot le  contraire n'est pas d\'emontr\'e dans \cite{P-Tr},
que $\Si_X$ ou $X$ aient des singularit\'es en codimension un. 
Pour cette raison il n'est pas exclu que le germe de tissu dont 
il est question \`a la fin de l'\'enonc\'e soit d\'efini sur un germe  
$(\Si'_X,x_0)$, lisse  de dimension $rn$ mais distinct 
du germe de $\Si_X$ en $x_0$.

\bk
Notre dernier \'enonc\'e fait le lien entre les deux pr\'ec\'edents  
en montrant que la propri\'et\'e 
pour un tissu de rang maximal d'\^etre isomorphe \`a un germe 
de tissu alg\'ebrique grassmannien ne d\'epend que de sa vari\'et\'e 
de Blaschke. Il utilise une notion dont nous esquissons 
maintenant la d\'efinition, voir le d\'ebut du Chapitre 3 pour une pr\'esentation 
pr\'ecise.

\sk
Soit $\cal{T}$ un $d$-tissu semi-extr\'emal de type $(r,n)$
sur un germe $(M,x_0)$.
Nous montrerons qu'il d\'efinit une {\em structure presque 
grassmannienne de type $(r,n)$} sur $(M,x_0)$,
une propri\'et\'e d\'ej\`a pressentie mais qui n'avait \'et\'e \'etablie auparavant 
que dans quelques cas.
Il s'agit d'un type de ${\rm G}$-structure model\'e sur la g\'eom\'etrie 
de la grassmannienne $\G_{r,n}$ des $(n-1)$-plans 
de $\P^{r+n-1}$. Pour guider l'intuition du lecteur, 
on peut traduire cette propri\'et\'e de la mani\`ere suivante.
En chaque point $x$ de $(M,x_0)$, le tissu $\cal{T}$
induit un $d$-tissu constant de m\^eme type sur l'espace tangent $T_xM$ et
celui-ci est localement isomorphe \`a un germe de tissu grassmannien.

Cette propri\'et\'e n'apporte rien si $r=1$ mais elle joue un r\^ole 
important si $r\geq 2$, ce que nous supposons 
maintenant. 
C'est un r\'esultat classique d'Akivis \cite{Ak},
qui appelle un $d$-tissu $\cal{T}$ avec $d\geq n+1$,
qui a la propri\'et\'e ci-dessus, un tissu presque grassmannien,
qu'un tel tissu est isomorphe 
\`a un germe de tissu grassmannien si et seulement si la structure 
presque grassmannienne qu'il d\'efinit est int\'egrable, c'est-\`a-dire localement
isomorphe \`a celle de la grassmannienne $\G_{r,n}$.

Si c'est le cas et s'il poss\`ede une relation ab\'elienne compl\`ete,
il est isomorphe \`a un germe 
de tissu alg\'ebrique grassmannien compte tenu du th\'eor\`eme d'Abel inverse.

D'autre part, nous avons montr\'e dans \cite{P-Tr} que, si $X$ est une 
vari\'et\'e de la classe $\cal{X}_{r+1,n}(q)$, la vari\'et\'e lisse $\Si_{X, {\rm adm}}$ 
des points admissibles de $\Si_X$ est naturellement munie d'une structure 
presque grassmannienne de type $(r,n)$. Dans le cas o\`u $X$ est la vari\'et\'e 
de Blaschke d'un tissu de rang maximal $\cal{T}$, cette structure co\"{\i}ncide 
avec celle que d\'efinit par transport le germe de tissu isomorphe \`a $\cal{T}$
donn\'e par le Th\'eor\`eme \ref{Th1B}.

\sk
Notre dernier \'enonc\'e est le suivant :
\bt
\label{Th1C}
Soit $\cal{T}$ un $d$-tissu de type $(r,n)$ et de rang maximal, avec $r\geq 2$ et $q(d)\geq n-1$,
et $X_{\cal{T}}$ sa vari\'et\'e de Blaschke. Le tissu $\cal{T}$
est isomorphe \`a un germe d'un tissu alg\'ebrique grassmannien 
si et seulement si la structure presque grassmannienne de $\Si_{X_{\cal{T}},{\rm adm}}$ est int\'egrable.
\et
Les trois th\'eor\`emes pr\'ec\'edents sont d\'emontr\'es dans le Chapitre 4, 
dans l'ordre suivant. Nous d\'emontrerons d'abord le Th\'eor\`eme \ref{Th1B}.
La d\'emonstration ne d\'epend que de r\'esultats relativement simples
de \cite{P-Tr} et d'un analogue du Th\'eor\`eme d'Abel inverse
pour les germes de tissus d'incidence pour une paire $(X,\Si_X)$. 
Le Th\'eor\`eme \ref{Th1C} sera alors une cons\'equence de la 
d\'efinition dans \cite{P-Tr} de la structure presque grassmannienne 
de $\Si_{X,{\rm adm}}$.
Finalement le Th\'eor\`eme \ref{Th1A} est une cons\'equence directe 
d'un r\'esutat essentiel de \cite{P-Tr}, que cette structure 
est toujours int\'egrable si $n=2$ ou si $q\neq 2n-3$.

\bk
Il reste \`a dire quelques mots sur un probl\`eme int\'eressant 
que nos \'enonc\'es laissent en suspens : \'etant donn\'es $r\geq 2$
et $n\geq 3$, existe-t-il un $d$-tissu de rang maximal, de type $(r,n)$
avec $q(d)=2n-3$, non isomorphe \`a un germe de tissu 
alg\'ebrique grassmannien ? 

Une condition n\'ecessaire est qu'il existe 
une vari\'et\'e $X$ de la classe $\cal{X}_{r+1,n}(2n-3)$
telle que la structure presque grassmannienne de $\Si_{X,{\rm adm}}$
ne soit pas int\'egrable.
Malheureusement nous savons peu de choses \`a ce sujet,
seulement qu'il en existe si $n$ est \'egal \`a trois ou \`a quatre,
quel que soit $r\geq 2$, ainsi que si $r=2$ et $n\in \{5,6\}$.

Dans un article \`a venir \cite{P-Tr2}, nous montrerons 
que, pour ces valeurs de $(r,n)$ et $q(d)=2n-3$,
il existe des $d$-tissus de rang maximal qui ne sont pas 
isomorphes \`a des germes  de tissus alg\'ebriques grassmanniens.

\ssct{Notes sur l'origine du probl\`eme}

La th\'eorie ab\'elienne des tissus n'a d'abord concern\'e que les tissus plans,
autrement dit le  cas $(r,n)=(1,2)$,  que nous excluons de nos 
\'enonc\'es. 

\sk
En 1932,  Blaschke et Howe obtiennent qu'un $d$-tissu plan
grassmannien sur $(\P^{2\star},x_0)$, {\em i.e.}
dont les feuilles sont des morceaux de droites de $\P^{2\star}$,
qui poss\`ede une relation ab\'elienne 
compl\`ete, est alg\'ebrique. 
C'est la version initiale du th\'eor\`eme d'Abel inverse exprim\'ee
en termes de tissus.

\sk
Les auteurs introduisent alors le probl\`eme de \og l'alg\'ebrisation \fg\,
des tissus plans de rang maximal. Dans ce qui suit, \og alg\'ebrisable \fg\, 
voudra dire isomorphe \`a un germe de tissu alg\'ebrique grassmannien.

\sk
La borne $\rho_{1,2}(d)=(d-1)(d-2)/2$ est obtenue par Bol.
Il est facile de v\'erifier que tous les $3$-tissus plans de rang maximal $1$
sont localement isomorphes. Ils sont donc alg\'ebrisables 
et m\^eme d'une infinit\'e de fa\c{c}ons puisqu'il existe une infinit\'e
de cubiques de $\P^2$ deux-\`a-deux non \'equivalentes.

\sk
La m\^eme ann\'ee, Howe r\'esout le probl\`eme pour les $4$-tissus 
de rang maximal $3$, en traduisant le th\'eor\`eme de Lie en termes de tissus.

\sk
En 1933 Blaschke compl\`ete la m\'ethode canonique pour r\'esoudre 
le probl\`eme dans le cas des $5$-tissus de rang maximal~$6$. 
Sa d\'emonstration est incorrecte et Bol publie 
en 1936 un contre-exemple, montrant que l'\'enonc\'e aussi est faux. 

\sk
Concernant les tissus plans, aucun progr\`es n'est \`a signaler juqu'en 2002--2003, 
quand Pirio et Robert construisent,
ind\'ependamment, de nouveaux contre-exemples. 
Le r\'esultat 
suivant a \'et\'e obtenu par Mar\'{\i}n, Pereira et Pirio \cite{M-Pe-P}.
\bt
Pour tout $d\geq 5$, il existe des $d$-tissus plans de rang maximal 
qui ne sont pas alg\'ebrisables.
\et
Nous renvoyons au livre Pereira-Pirio \cite{Pe-P} pour une 
pr\'esentation de la th\'eorie des tissus 
plans et en particulier des r\'esultats r\'ecents.
On y trouvera aussi plus d'informations qu'ici sur la gen\`ese de la 
th\'eorie, en particulier du th\'eor\`eme de Blaschke-Howe
mais, sur cette gen\`ese, le livre Blaschke-Bol \cite{BB} 
reste bien s\^ur la r\'ef\'erence fondamentale.

\bk
Venons-en aux tissus de codimension un avec $n\geq 3$.
En 1934 Bol obtient, par la m\'ethode canonique, le r\'esultat majeur 
qu'un $d$-tissu $\cal{T}$, de type $(1,3)$
et de rang maximal,  est alg\'ebrisable si $d\geq 6$. 
L'argument de Bol pour d\'emontrer que ce que nous appelons la vari\'et\'e de Blaschke du tissu $\cal{T}$
est une surface est en partie g\'eom\'etrique, bas\'e sur une analogie avec une g\'eom\'etrie de Weyl.

\sk
Chern a d\'etermin\'e en 1935 la borne $\rho_{1,n}(d)$ pour le rang d'un
tissu de codimension $1$. Plus de quarante ans plus tard, il entreprend 
avec Griffiths d'\'etendre le r\'esultat de Bol en toute dimension $n\geq 3$.

Les deux premi\`eres parties de Chern-Griffiths \cite{C-Gr1} restent une 
r\'ef\'erence essentielle pour qui veut comprendre  
les liens entre le  probl\`eme de l'alg\'ebrisation des tissus et la g\'eom\'etrie 
alg\'ebrique. 
Dans la troisi\`eme partie, de plus de quarante pages, les auteurs traitent le noyau dur 
du probl\`eme : montrer que la vari\'et\'e de Blaschke est une surface.
Comme chez Bol, l'id\'ee \`a la base des calculs qu'ils font est de nature 
g\'eom\'etrique, la g\'eom\'etrie des chemins se substituant \`a la
g\'eom\'etrie de Weyl. Ils  reconna\^{\i}tront trois ans plus tard dans 
\cite{C-Gr2}  une erreur dans cette partie de l'article, qui les force
\`a affaiblir de fa\c{c}on consid\'erable leur \'enonc\'e initial.
 
\sk
Le probl\`eme est finalement r\'esolu par Tr\'epreau \cite{Tr},
qui obtient le r\'esultat suivant.
\bt
Un $d$-tissu semi-extr\'emal de type $(1,n)$ avec $n\geq 3$ est alg\'ebrisable.
\et
L'argument de \cite{Tr}, pour montrer que la vari\'et\'e de Blaschke est une surface,
est analytique. Il s'agit d'un calcul, laborieux mais direct, du rang d'un 
param\'etrage local de cette vari\'et\'e.
La m\^eme id\'ee appara\^{\i}t d\'ej\`a dans l'article 
de H\'enaut \cite{He} \`a propos des tissus de type $(r,2)$,
un cas pour lequel le calcul du rang de ce param\'etrage est plus facile.

\bk
En 1976 para\^{\i}t l'article de Griffiths,  {\em Variations on a theorem of Abel}~\cite{Gr},
qui donne des versions g\'en\'erales du th\'eor\`eme d'Abel
et du th\'eor\`eme d'Abel inverse de Darboux. En 1978, la 
m\^eme ann\'ee que~\cite{C-Gr1}, para\^{\i}t 
l'article de Chern et Griffiths \cite{C-Gr3} sur le rang des tissus 
de type $(r,n)$. Outre la borne sur le rang, ils y d\'ecrivent la g\'eom\'etrie
des normales d'un tissu de rang maximal de codimension~$2$.
La d\'ecouverte d'une lacune dans \cite{C-Gr1} est peut-\^etre l'une des raisons 
pour lesquelles ils n'ont pas poursuivi l'\'etude de ce probl\`eme.

\sct{Tissus constants semi-extr\'emaux}

Dans ce chapitre, nous ne consid\'erons que des {\em tissus constants}, 
de type $(r,n)$ sur un espace vectoriel $V$ de dimension $rn$,
c'est-\`a-dire des tissus dont les \'el\'ements sont des 
feuilletages de $V$ en sous-espaces affines parall\`eles
de codimension $r$.

Nous verrons au d\'ebut du Chapitre 3 que,  pour ce qui est des questions que nous 
abordons ici, on ne perd rien en g\'en\'eralit\'e \`a supposer les tissus constants. 
Et l'on gagne beaucoup quant \`a la clart\'e des notations.

\sk
Nous rappellerons d'abord la d\'emonstration du Th\'eor\`eme \ref{borne},
la borne de Chern-Griffiths, que nous \'enoncerons sous une forme pr\'ecis\'ee.

Nous caract\'eriserons ensuite, c'est le r\'esultat important de ce chapitre,
la propri\'et\'e pour un tissu constant d'\^etre semi-extr\'emal par la g\'eom\'etrie de son syst\`eme 
de normales, voir la Proposition \ref{tensor}. 

\ssct{Les  bornes de Chern-Griffiths}

Soit $\cal{T}=\{\cal{F}_1,\ldots,\cal{F}_d\}$ 
un tissu constant de type $(r,n)$ sur l'espace vectoriel $V$.
Les  feuilles de $\cal{F}_j$
sont les fibres d'une application lin\'eaire de rang $r$
$$
{\k}_j = (u_{j,1},\ldots,u_{j,r}): \, V\rightarrow \C^r
$$
et $\,\Omega_j = du_{j,1}\wedge \cdots \wedge du_{j,r}$
est une normale g\'en\'eratrice de $\cal{F}_j$.
Comme la normale $\O_j$ est ferm\'ee, une normale 
$\phi_j = c_j\Omega_j$ de $\cal{F}_j$  sur $(V,0)$ est 
ferm\'ee si et seulement si $dc_j \wedge \O_j=0$,
donc si et seulement si la fonction $c_j$ est de la forme 
$c_j=f_j({\k}_j)$, o\`u $f_j$ est une fonction sur $(\C^r,0)$.
De plus, le germe $c_j$ est nul si et seulement si le germe
$f_j$ est nul.

\sk
Avec ces notations, une relation ab\'elienne de $\cal{T}$ sur $(V,0)$ est 
un $d$-uplet de normales ferm\'ees $f_j(\k_j)\O_j$ dont la somme est nulle.
Soit $f_j=\sum_{h=0}^{+\infty} p_{j,h}$ 
la d\'ecompositon 
de $f_j$ en s\'erie de polyn\^omes homog\`enes,  $p_{j,h}$
\'etant de degr\'e $h$. Comme les formes $\O_j$ sont constantes, il est clair
que le $d$-uplet des $p_{j,h}(\k_j)\O_j$
est une relation ab\'elienne du tissu $\cal{T}$, homog\`ene de degr\'e~$h$.

\sk
Suivant \cite{C-Gr3} nous allons d\'emontrer la forme 
pr\'ecis\'ee suivante de (\ref{borne}).
\bpr
\label{CG-cste}
Soit $\cal{T}$ un $d$-tissu constant de type $(r,n)$ sur l'espace vectoriel $V$.
Pour tout $h\in \N$, la dimension de l'espace $R(h)$ des relations ab\'eliennes de $\cal{T}$
qui sont homog\`enes de degr\'e $h$ v\'erifie 
\beq
\label{borne-h}
\dim R(h) \leq {r-1+h \choose r-1} \times \max\,\left( d-(r+h)(n-1)-1,0\right).
\eeq
\epr
Nous montrerons dans la prochaine section qu'il existe des tissus constants 
pour lesquels les bornes (\ref{borne-h}) sont atteintes quel que soit $h$.

\sk
Soit $E_r(h)$ l'espace, de dimension ${r-1+h\choose r-1}$,
des polyn\^omes 
homog\`enes de degr\'e $h$ sur $\C^r$ et $N(h)$ l'espace engendr\'e par les normales ferm\'ees 
homog\`enes de degr\'e $h$ du tissu $\cal{T}$. L'application lin\'eaire 
$E_r(h)^d \rightarrow N(h)$, d\'efinie par 
$(c_1,\ldots,c_d)  \mapsto \; \sum_{j = 1}^d c_j({\k}_j) \, \Omega_j$,
est surjective et son noyau est isomorphe \`a l'espace $R(h)$.
On a donc 
$$
\dim R(h) + \dim N(h) = d\dim E_r(h).
$$

\sk
Pour d\'emontrer (\ref{borne-h}), il suffit de montrer qu'on a :
\beq
\label{borne-hbis}
d = (r+h)(n-1)+1 \; \Rightarrow \; R(h)=0,
\eeq
ce qui donne alors $\dim N(h) = ((r+h)(n-1)+1)\dim E_r(h)$.
En effet on en d\'eduit qu'on a encore $R(h)=0$ si 
$d\leq (r+h)(n-1)+1$ et, comme la dimension minimale de $N(h)$
ne peut que cro\^{\i}tre avec $d$, qu'on a :
$$
d > (r+h)(n-1) +1 \; \Rightarrow \; \dim R(h) \leq (d-(r+h)(n-1)-1) \dim E_r(h).
$$
\bpf[D\'emonstration de (\ref{borne-hbis}) pour $h=0$]
Soit $(c_1\O_1,\ldots,c_d\O_d)$
une relation ab\'elienne \`a coefficients constants.
Par sym\'etrie, il suffit de montrer que $c_1$ est nul.

\sk
Pour le voir, on multiplie la relation $\sum_{j=1}^d c_j\O_j=0$ par une $(d-1)$-forme $K$
telle que $\O_1\wedge K\neq 0$ et $\O_j\wedge K=0$ si $j\geq 2$, ce qui donne 
$c_1(\O_1\wedge K) = (\sum_{j=1}^d c_j\O_j) \wedge K = 0$
et le r\'esultat. 

\bk
Reste \`a v\'erifier qu'une telle forme $K$ existe. (C'est consid\'er\'e 
comme \'evident dans \cite{C-Gr3}.) 
Pour les besoins de la r\'ecurrence, nous supposons seulement 
$d\geq n$ et $d\leq r(n-1)+1$. L'existence de $K$ est \'evidente 
si $d=n$ puisque la condition (PG) donne alors $\O_1\wedge \cdots \wedge \O_n\neq 0$.

Notons $F_j=\text{Ker}\, \kappa_j \subset V$ la direction du feuilletage $\cal{F}_j$
et $F_j^\perp \subset V^\star$ son orthogonal. \'Etant donn\'e
$v=(v_2,\ldots,v_d)\in \prod_{j=2}^d F_j^\perp$,
posons 
$$
K(v) =\wedge_{j=2}^d \, dv_j, \qquad K_k(v) = \wedge_{j=2, j\neq k}^d \, dv_j, \qquad k=2,\ldots,d. 
$$
On a toujours $\O_j\wedge K(v)=0$ si $j\geq 2$. D'autre part, 
pour un tissu donn\'e,  $\O_1\wedge K(v)$ est ou bien toujours nul,
ou bien non nul pour tout $v$ dans un ouvert dense de $\prod_{j=2}^d F_j^\perp$.
Pour montrer que c'est la seconde propri\'et\'e qui est vraie,
nous raisonnons par l'absurde et par r\'ecurrence, en supposant $d>n$ et
que la seconde propri\'et\'e est vraie pour tous les sous-tissus 
d'ordre $d-1$ de $\cal{T}$. Plus pr\'ecis\'ement nous pouvons supposer
\beq
\label{banal}
\O_1\wedge K(v) = 0, \;\;\; \O_1\wedge K_k(v) \neq 0 \;\;\; (k=2,\ldots,d),
\eeq
pour tout $v\in \prod_{j=2}^d F_j^\perp$ voisin d'un certain $v^0$ fix\'e.
Soit $H$ le sous-espace de $V^\star$, de dimension $r+d-2 < rn$, engendr\'e par $F_1^\perp$
et les formes $v_2^0,\ldots,v_{d-1}^0$. 

\sk
En appliquant (\ref{banal}) \`a $v^0$, nous obtenons que $v_d^0$
appartient \`a $H$. En \'ecrivant (\ref{banal}) pour $v=v^0$
puis, pour $k\in \{2,\ldots,n\}$ donn\'e,
en faisant varier $v_k$ au voisinage de $v^0_k$, en particulier
$v_k$ d\'ecrit une base de $F_k^\perp$, 
nous obtenons que $F_k^\perp$ est contenu dans $H$.

Ainsi $H$ contient $F_1^\perp,\ldots,F_n^\perp$. C'est une contradiction 
car, compte tenu de la condition de position g\'en\'erale (PG),  ces espaces sont en 
somme directe et donc engendrent un espace de dimension 
$rn>\dim H$.
\epf

\sk
\bpf[D\'emonstration de (\ref{borne-hbis}) pour $h\geq 1$]
Nous supposons que (\ref{borne-hbis}) est 
vrai \`a l'ordre $h-1$ et nous le d\'emontrons \`a l'ordre $h$. 
Par sym\'etrie, il suffit de montrer que la $n$-i\`eme composante 
d'un \'el\'ement de $R(h)$ est nulle. 

\sk
Le sous-tissu $\cal{T}'=\cal{T}\bck \{\cal{F}_1,\dots,\cal{F}_{n-1}\}$
est d'ordre $d-(n-1)$ donc, par hypoth\`ese de r\'ecurrence, 
l'espace $R'(h-1)$
de ses relations homog\`enes de degr\'e $h-1$ est r\'eduit \`a $0$.

Soit $X\in \cap_{j=1}^{n-1}F_j$, vu comme un champ de vecteurs constant sur $V$.
Notons que 
$$
X\cdot c_j(\k_j) = \sum_{a=1}^r (X\cdot u_{j,a}) \fr{\pl c_j}{\pl t_a}(\k_j),
$$
o\`u les coefficients sont des scalaires. Il en r\'esulte 
que $(X\cdot c_j(\k_j))\O_j$ est une normale ferm\'ee de $\cal{F}_j$,
nulle si $j=1,\ldots,n-1$. En appliquant un tel vecteur \`a un
\'el\'ement de $R(h)$, on obtient donc en fait un \'el\'ement 
de $R'(h-1)$, donc une relation triviale.

En particulier, la $n$-i\`eme composante de la relation initiale 
est annul\'ee \`a la fois par les vecteurs $X\in F_n$
et par les \'el\'ements du suppl\'ementaire $\cap_{j=1}^{n-1} F_j$
de $F_n$. Cette composante est constante et homog\`ene de degr\'e $h\geq 1$,
donc nulle.
\epf
Finissons par une remarque. Soit $\cal{T}'$ un sous-tissu d'ordre $d'=d-1$
du $d$-tissu $\cal{T}$, avec $d-1\geq (r+h)(n-1)+1$.
Notons $R(h)$, $R'(h)$,
$N(h)$ et  $N'(h)$ les espaces associ\'es de relations 
et de normales ferm\'ees homog\`enes de degr\'e $h$.
On sait que  $\dim R(h) + \dim N(h) = d \dim E_r(h)$
et que $\dim R'(h) + \dim N'(h) = (d-1)\dim E_r(h)$.

Si la dimension de $R(h)$ est maximale, comme $N'(h)$ est un sous-espace 
de $N(h)$, on obtient $\dim R'(h) \geq \dim R(h) - \dim E_r(h)$
et donc que $R'(h)$ est de dimension maximale.

En appliquant cette remarque \`a $h\in \{0,1\}$ et 
en it\'erant, on obtient le lemme suivant, que nous utiliserons.
\ble
\label{sous-tissu}
Un sous tissu $\cal{T}'$,  d'ordre $d'$ avec $q(d')\geq n-1$, d'un $d$-tissu
{\em constant} semi-extr\'emal est aussi semi-extr\'emal.
\ele

\ssct{Exemples}

Nous continuons avec les notations pr\'ec\'edentes.

Notons $A=\{1,\ldots,r\}\times \{1,\ldots,n\}$.
Soit $m = (m_{a,\a})_{(a,\a)\in A}$ une base de $V^\star$.
Cette base \'etant fix\'ee, nous associons \`a 
tout point ${p} = [\xi_1 : \cdots : \xi_n]$ de $\P^{n-1}$, 
le feuilletage $\cal{F}(p)$, constant de type $(r,n)$,
dont la direction est donn\'ee par le syst\`eme d'\'equations 
$$
\sum_{\a=1}^n \xi_\a \, m_{a,\a} =0, \qquad  a=1,\ldots,r.
$$
Ces feuilletages n'ont rien de particulier si $r=1$.
Nous supposons maintenant $r\geq 2$.
Des feuilletages $\cal{F}({p}_j)$ v\'erifient la condition 
(PG) si et seulement si les points ${p_j}$ 
sont en position g\'en\'erale dans $\P^{n-1}$.
Il suffit de le v\'erifier pour $d=n$ ; les deux propri\'et\'es se traduisent
alors par le fait que la matrice $n\times n$ form\'ee \`a partir 
des syst\`emes de composantes homog\`enes des $p_j$ est inversible.

\sk
Le lemme suivant est standard.
\ble
Une base $(m_{a,\a})_{(a,\a)\in A}$ de $V^\star$ \'etant donn\'ee,
si $d$ points deux-\`a-deux distincts $p_j$
appartiennent \`a une m\^eme 
courbe rationnelle normale de degr\'e $n-1$ de $\P^{n-1}$,
les bornes (\ref{borne-h}) sont atteintes pour tout $h$
par le tissu $\{\cal{F}(p_1),\ldots,\cal{F}(p_d)\}$.
Il est donc de rang maximal.
\ele
\bpf
Une courbe rationnelle normale de degr\'e $n-1$ dans $\P^{n-1}$
admet un param\'etrage de la forme $p(t) = [\phi_1(t) : \cdots : \phi_n(t)]$,
o\`u les $\phi_\a(t)$ forment une base de l'espace des polyn\^omes de degr\'e 
$\leq n-1$. Les points $p_j = p(t_j)$, on peut supposer
$t_j\in \C$, sont en position g\'en\'erale dans $\P^{n-1}$ 
d\`es que les $t_j$ sont deux-\`a-deux distincts.

\sk
Pour toute valeur de $t\in \C$, notons $l_a(t)= \sum_{\a=1}^n \phi_\a(t)\, m_{a,\a}$.
En d\'eveloppant $\O(t)= \wedge_{a=1}^r dl_a(t)$ suivant les puissances de $t$, on obtient
$$
\O(t) = \sum_{\rho=0}^{r(n-1)} t^\rho \, K_\rho,
$$
o\`u les $r$-formes constantes $K_\rho$ ne d\'ependent pas de $t$.
Ainsi  les normales g\'en\'eratrices constantes $\O(t)$
des feuilletages  $\cal{F}(p(t))$ engendrent, quand $t$ d\'ecrit $\C$,
un espace de dimension $\leq r(n-1)+1$.

\sk
Plus g\'en\'eralement, \'etant donn\'e $h\in \N$ et $c\in E_r(h)$, 
la fonction $c(l_1(t),\ldots,l_r(t))$
est polynomiale en $t$ de degr\'e $\leq h(n-1)$ donc 
$$
c(l_1(t),\ldots,l_r(t))\, \O(t) = \sum_{\rho=0}^{(r+h)(n-1)} t^\rho \, K_{c,\rho},
$$
o\`u les $K_{c,\rho}$ sont des $r$-formes homog\`enes de degr\'e $h$
sur $V$, ind\'ependantes de $t$ et d\'ependant lin\'eairement de $c\in E_r(h)$.
Quand $t$ d\'ecrit $\C$ et que $c$ parcourt une base de $E_r(h)$,
elles engendrent un espace 
de dimension $\leq \dim E_r(h) \times ((r+h)(n-1)+1)$.

\sk
Ceci d\'emontre le lemme. En effet, si le $d$-tissu $\cal{T}$ 
est comme dans l'\'enonc\'e, on obtient que l'espace $N(h)$
engendr\'e par ses normales ferm\'ees homog\`enes de degr\'e $h$ v\'erifie 
$$
\dim N(h) \leq \dim E_r(h) \times \min (d,(r+h)(n-1)+1),
$$
ce qui \'equivaut au fait que le premier membre de (\ref{borne-h}) 
est aussi plus grand que le second.
\epf
Le lemme pr\'ec\'edent a l'int\'er\^et de montrer simplement que la borne 
$\rho_{r,n}(d)$ pour le rang d'un tissu $\cal{T}$
est atteinte et, dans le cas d'un tissu constant,
qu'elle est atteinte si et seulement si les  in\'egalit\'es 
(\ref{borne-h}) sont des \'egalit\'es pour tout $h$.

Comme cons\'equence, si $q(d)\geq n-1$, un $d$-tissu constant
est semi-extr\'emal si et seulement si 
ses espaces $R(0)$ et $R(1)$ sont de dimension 
maximale, soit $d-r(n-1)-1$ et $r[d-(r+1)(n-1)-1]$
respectivement, et c'est le cas s'il est de rang maximal.

\sk
Le lemme pr\'ec\'edent admet la r\'eciproque suivante,
d\'ej\`a remarqu\'ee dans Chern-Griffiths \cite{C-Gr2}
et  Little \cite{Li}.
\ble
\label{castel-2}
On suppose $d\geq r(n-1)+ 1$ si $r\geq 3$,
$d\geq 2n+1$ si $r=2$.
Si le syst\`eme des normales constantes du 
tissu $\{\cal{F}({p}_1), \ldots , \cal{F}({p}_d)\}$ 
est de rang minimal $r(n-1)+1$, les points ${p}_j$ 
sont situ\'es sur une m\^eme courbe rationnelle normale de
degr\'e $n-1$ de $\P^{n-1}$.
\ele
\bpf
Si ${p} = [\xi_1 : \cdots : \xi_n]\in \P^{n-1}$, 
$$
\Omega({p}) = \bigwedge_{a=1}^r 
\left( \sum_{\alpha=1}^n  \xi_\alpha \, dm_{a,\alpha} \right)
 =  
\sum_{\alpha_1,\ldots,\alpha_r=1}^n
\xi_{\alpha_1}\ldots \xi_{\alpha_r} \, 
dm_{1, \alpha_1} \wedge \cdots \wedge dm_{r, \alpha_r}
$$
est une normale constante g\'en\'eratrice de $\cal{F}({p})$.
Tous les mon\^omes homog\`enes de degr\'e $r$ en $(\xi_1,\ldots,\xi_n)$
sont des composantes de $\O(p)$ dans la base canonique de $\wedge^r V^\star$,
avec des r\'ep\'etitions, les autres composantes \'etant nulles.
On en d\'eduit que le rang du syst\`eme des
normales des feuilletages $\cal{F}(p_j)$
est \'egal \`a la dimension de l'espace engendr\'e par les points 
$v_r({p_1}),\ldots v_r({p_d})$,  o\`u $v_r$ est un plongement
de Veronese d'ordre $r$ de $\P^{n-1}$, 
c'est-\`a-dire associ\'e au syst\`eme lin\'eaire $|\cal{O}_{\P^{n-1}}(r)|$.

\sk
C'est un r\'esultat classique de Castelnuovo, voir par exemple Harris~\cite{Ha},
que si $d\geq r(n-1)+1$ (respectivement $d\geq 2n+1$ si $r=2$),
la dimension de l'espace engendr\'e par les points $v_r(p_j)$
est de dimension minimale $r(n-1)+1$ si et seulement si les points $p_j$,
suppos\'es \^etre en position g\'en\'erale dans $\P^{n-1}$,
appartiennent \`a une m\^eme courbe rationnelle normale de degr\'e $n-1$.
\epf

\ssct{La g\'eom\'etrie du syst\`eme des normales}

Nous avons le r\'esultat suivant.
\bpr
\label{tensor}
Soit $\cal{T}$ un $d$-tissu constant semi-extr\'emal de type $(r,n)$
sur l'espace $V$ avec $r\geq 2$. Il existe une base $(m_{a,\a})_{(a,\a)\in A}$ de $V^\star$ 
et des points ${p}_1,\ldots,{p}_d$, en position g\'en\'erale dans $\P^{n-1}$,
tels que les composantes de $\cal{T}$ soient les feuilletages 
$\cal{F}({p}_j)$. De plus les points $p_j$ sont situ\'es sur une m\^eme 
courbe rationnelle normale de degr\'e $n-1$ dans $\P^{n-1}$.
\epr
Chern et Griffiths cherchent dans \cite{C-Gr2}
\`a \'etablir un r\'esultat de cette nature. Ils l'obtiennent
pour $r=2$, mais en renfor\c{c}ant la condition
(PG) si $n\geq 3$. Ils n'utilisent que le fait que le syst\`eme 
des normales {\em constantes} est de rang minimal $2(n-1)+1$.
Little cherche dans \cite{Li} \`a g\'en\'eraliser l'\'etude 
de \cite{C-Gr2} pour $r\geq 3$ mais sa d\'emonstration semble reposer 
sur un cercle vicieux, \`a la fin de la page 26
ou au d\'ebut de la page 27.
\bpf
La seconde partie de l'\'enonc\'e est une cons\'equence de la premi\`ere 
et du Lemme \ref{castel-2}.

Notons $\cal{F}_j$ les composantes du tissu $\cal{T}$,
$F_j$ la direction et $\O_j$ une normale g\'en\'eratrice constante de $\cal{F}_j$.

\`A l'oppos\'e des travaux que nous venons de mentionner,
c'est \`a partir de l'espace des relations ab\'eliennes {\em homog\`enes de degr\'e~$1$}
du tissu que nous allons construire la base adapt\'ee $(m_{a,\a})_{(a,\a)\in A}$.

\sk
\nk
{\em Nous supposons d'abord $q(d)=n-1$.} 

\sk
Alors $d = (r+1)(n-1) + 2$ et par hypoth\`ese l'espace $R(1)$ des relations 
homog\`enes de degr\'e $1$ est de dimension maximale 
$r(d-(r+1)(n-1)-1)$ c'est-\`a-dire de dimension $r$. Soit  
$$
(u_{a,1} \, \O_1,\ldots,u_{a,d} \,\O_d), \qquad a=1,\ldots,r,
$$
une base de $R(1)$. Nous affirmons que $F_j^\perp = \lan u_{1,j},\ldots,u_{r,j} \ran$
pour tout $j$.  Sinon
une combinaison lin\'eaire non triviale 
convenable des relations ci-dessus aurait sa $j$-i\`eme 
composante nulle. Elle induirait une relation ab\'elienne non triviale
du $(d-1)$-tissu $\cal{T}\bck\{\cal{F}_j\}$, ce qui est impossible
d'apr\`es la majoration (\ref{borne-h}) appliqu\'ee \`a $d-1$. 

\sk
Comme $V^\star = \oplus_{j=1}^n F_j^\perp$, nous d\'efinissons une base de $V^\star$
en posant
$$
m_{a,\a} = u_{a,\a}, \qquad a=1,\ldots,r, \;\; \a=1,\ldots,n.
$$
Pour ce choix de base et par construction, pour $j=1,\ldots,n$ le feuilletage $\cal{F}_j$
est le feuilletage $\cal{F}(p_j)$ associ\'e au point $p_j\in \P^{n-1}$ 
dont les coordonn\'ees homog\`enes sont toutes nulles sauf la $j$-i\`eme.

D'autre part, les normales $\O_j$ engendrent un espace de dimension $r(n-1)+1=d-n$
et sont en position g\'en\'erale dans cet espace. Nous 
choisissons $\O_{n+1},\ldots,\O_d$ comme base de cet espace et 
d\'ecomposons 
$$
\Omega_\a = \sum_{j=n+1}^d \xi_{\a, j} \, \Omega_j,
\qquad 
\a= 1,\ldots,n.
$$
Les relations $\sum_{j=1}^d u_{a,j}\Omega_j = 0$ s'\'ecrivent alors 
$$
\sum_{j=n+1}^d \left( u_{a,j} + \sum_{\a=1}^n \xi_{\a,j} u_{a,\a}   \right) \Omega_j = 0.
$$
Si $j\geq n+1$, le feuilletage $\cal{F}_j$ est donc d\'efini par le syst\`eme
$$
\sum_{\a=1}^n \xi_{\a,j} m_{a,\a}=0, \qquad a=1,\ldots,r.
$$
C'est le feuilletage $\cal{F}(p_j)$ associ\'e au point $p_j=[\xi_{1,j} : \cdots : \xi_{n,j}]$
de $\P^{n-1}$ pour la base des $m_{a,\a}$. Ceci d\'emontre la proposition 
dans le cas particulier o\`u $q(d)=n-1$.

\sk
\nk
{\em Nous supposons maintenant $q(d)\geq n-1$.} 

\sk
Pour en d\'eduire le r\'esultat dans le cas g\'en\'eral, le plus simple 
est d'appliquer le Lemme \ref{Akivis}, d\^u \`a Akivis \cite{Ak}, que nous 
d\'emontrerons au d\'ebut du Chapitre 3 : un $(n+1)$-tissu
{\em d\'etermine} sur $V$ une structure presque grassmannienne constante.

Nous venons de montrer 
que si $q(d)=n-1$, un $d$-tissu semi-extr\'emal est presque grassmannien,
c'est-\`a-dire que les structures d\'etermin\'ees 
par ses sous-tissus d'ordre $(n+1)$ co\"{\i}ncident.

Dans le cas g\'en\'eral, comme tous les sous-tissus de $\cal{T}$,  d'ordre $d'$
avec $q(d')\geq n-1$, sont aussi semi-extr\'emaux d'apr\`es le Lemme \ref{sous-tissu},
il suffit d'appliquer le cas particulier d\'ej\`a r\'esolu 
\`a tous les sous-tissus 
d'ordre $(r+1)(n-1)+2$ qui contiennent $\cal{F}_1,\ldots,\cal{F}_{n+1}$
pour obtenir encore le r\'esultat.
\epf

\sct{Tissus semi-extr\'emaux et la m\'ethode canonique}

Nous revenons au cas g\'en\'eral. Nous rappelons d'abord la notion
de structure presque grassmannienne, nous dirons aussi ${\rm G}_{r,n}$-structure,
et sa relation avec la th\'eorie ab\'elienne des tissus,
suivant Akivis \cite{Ak}.

\sk
\'Etant donn\'e un tissu $\cal{T}$ sur un germe $(M,x_0)$,
nous d\'efinissons son tissu tangent $\cal{T}_x$ en un point $x$.
C'est un tissu constant sur l'espace tangent $T_xM$, dont le rang est
au moins \'egal \`a celui de $\cal{T}$. 
Ce fait nous permet d'appliquer les r\'esultats obtenus 
pour les tissus constants dans le chapitre pr\'ec\'edent,
d'abord bien s\^ur d'obtenir la borne de Chern-Griffiths 
dans le cas g\'en\'eral. 

\sk
La principale 
application concerne les tissus semi-extr\'emaux.
Il r\'esulte de la Proposition \ref{tensor}
qu'un tel tissu d\'etermine une structure 
presque grasmannienne sur le germe $(M,x_0)$
o\`u  il est d\'efini, et m\^eme un peu plus.
Ces propri\'et\'es sont suffisantes pour mettre en \oe uvre 
la m\'ethode canonique de fa\c{c}on efficace.
Elles sont d'ailleurs inutiles dans le cas 
particulier o\`u $q(d)=n-1$, dont nous rappelons 
la d\'emonstration historique.

\sk
La fin du chapitre est plus technique. Il s'agira de calculer 
le rang d'une application qui param\`etre le germe de vari\'et\'e
d\'ecrit par les courbes $\s(x)$ dont il est question dans 
la D\'efinition \ref{var-X} de la vari\'et\'e de Blaschke.
Le r\'esultat cl\'e, que cette application est de rang 
$r+1$, nous permettra dans le prochain chapitre de montrer
que la vari\'et\'e de Blaschke $X_\cal{T}$ du tissu 
est elle-m\^eme de dimension $r+1$.

\ssct{Structures et tissus presque grassmanniens}

Nous notons \`a nouveau $A=\{1,\ldots,r\}\times \{1,\ldots,n\}$
et nous indexons les \'el\'ements d'une base de $\C^{rn}$
par les paires $(a,\a)\in A$. 
Les entr\'ees d'une  matrice $g\in {\rm GL}(\C^{rn})$ s'\'ecrivent alors 
$\,g_{a\a,b\b}\,$ avec $(a,\a),(b,\b)\in A$.

Soit ${\rm G}_{r,n}$ le sous-groupe des matrices $g\in {\rm GL}(\C^{rn})$ de la forme   
$$
g_{a \a,b \b} = C_{ab}\, A_{\a\b}, \;\;\; (a,\a), \, (b,\b)\in A.
$$ 
Une {\em structure presque grassmannienne de type $(r,n)$} ou {\em ${\rm G}_{r,n}$-structure} 
sur le germe 
lisse $(M,x_0)$ de dimension $rn$ est une ${\rm G}$-structure sur $(M,x_0)$
pour le groupe ${\rm G}={\rm G}_{r,n}$.

\sk
Parmi les diff\'erentes fa\c{c}ons de se donner une telle 
structure, la plus commode ici est de la d\'efinir 
par une famille maximale 
de bases du module des $1$-formes sur $(M,x_0)$
telle que les formules de passage d'une base 
$(m_{a,\a})_{(a,\a)\in A}$ de la famille \`a une autre $(n_{a,\a})_{(a,\a)\in A}$
sont de la forme 
\beq
\label{passage}
n_{a,\a} = \sum_{b=1}^r \sum_{\b =1}^n C_{ab} A_{\a \b} \, m_{b,\b},
\eeq
o\`u les coefficients $C_{ab}$ et $A_{\a\b}$ sont des fonctions sur $(M,x_0)$.
Nous dirons qu'une base qui appartient \`a la famille d\'efinit la structure.

Nous utiliserons les notions g\'eom\'etriques suivantes qu'une ${\rm G}_{r,n}$-structure permet de 
d\'efinir.
\bd
Soit $(M,x_0)$ un germe lisse muni de la ${\rm G}_{r,n}$-struc\-ture
d\'efinie par une base $(m_{a,\a})_{(a,\a)\in A}$.
Une {\em sous-vari\'et\'e distingu\'ee} de $(M,x_0)$
est une sous-vari\'et\'e $N$ de codimension~$r$, dont l'espace 
tangent $T_xN$ en un point est donn\'e par des \'equations de la forme 
\beq
\label{systeme-xi}
\sum_{\a=1}^n \xi_\a(x) m_{a,\a}(x) = 0, \qquad a=1,\ldots,r.
\eeq
\ed
\bd
Un feuilletage sur $(M,x_0)$ 
est {\em presque grassmannien} si ses feuilles sont des sous-vari\'et\'es distingu\'ees
de $(M,x_0)$. Un tissu est presque grassmannien si les feuilletages qui le composent
le sont.
\ed
Il est tr\`es facile de v\'erifier que ces notions  
ne d\'ependent  que de la structure et pas de la base d\'efinissante choisie.
D'autre part, la structure est d\'etermin\'ee 
par ses sous-vari\'et\'es distingu\'ees.
L'\'enonc\'e plus pr\'ecis suivant  est d\^u \`a Akivis \cite{Ak}.
Nous l'avons utilis\'e \`a la fin du chapitre pr\'ec\'edent.
\ble 
\label{Akivis}
Soit $\{\cal{F}_1,\ldots,\cal{F}_{n+1}\}$ un $(n+1)$-tissu  
de type $(r,n)$ sur un germe $(M,x_0)$.
Il existe une et une seule ${\rm G}_{r,n}$-structure sur $(M,x_0)$ 
pour laquelle ce tissu est presque grassmannien.
\ele
\bpf
Chaque feuilletage $\cal{F}_j$ est d\'efini par un syst\`eme
$$
\o_{a,j} = 0, \qquad a=1,\ldots,r.
$$
Compte tenu de la condition (PG), les $\o_{a,\a}$ avec $\a\neq n+1$
forment une base du module des $1$-formes sur $(M,x_0)$.
Pour chaque $a$, on d\'ecompose $\o_{a,n+1}$ dans cette base 
$$
\o_{a,n+1} = \sum_{\a=1}^n m_{a,\a}, \qquad m_{a,\a} = \sum_{b=1}^r C_{a b \a} \o_{b,\a}.
$$
La condition (PG) implique que les $m_{a,\a}$ forment une base.¨
Par construction, les feuilletages $\cal{F}_j$
sont d\'efinis respectivement par les syst\`emes 
$$
m_{1,j}=\cdots = m_{r,j} = 0 \;\; (j=1,\ldots,n) , \;\;\; \sum_{\a=1}^n m_{1,\a}= \cdots = \sum_{\a=1}^n m_{r,\a}= 0.
$$
En particulier ils sont presque grassmanniens pour la ${\rm G}_{r,n}$-structure 
d\'efinie par la base $(m_{a,\a})_{(a,\a)\in A}$.

\sk
Supposons qu'ils sont aussi presque grassmanniens 
pour la structure d\'efinie par une autre base $(n_{a,\a})_{(a,\a)\in A}$.
Par d\'efinition, chaque feuilletage $\cal{F}_j$ est d\'efini 
dans la nouvelle base 
par un syst\`eme analogue \`a (\ref{systeme-xi}) dont 
les coefficients forment un syst\`eme de coordonn\'ees
homog\`enes de $p_j(x)$ pour une certaine fonction $p_j: (M,x_0)\rightarrow \P^{n-1}$.
Un changement de base 
$$
n'_{a,\a} = \sum_{\b=1}^n A_{\a\b} n_{a,\b},
$$
ne modifie pas la structure et permet de supposer 
que $p_{n+1}=[1:\cdots:1]$ et que, si $j\leq n$,  les coordonn\'ees
homog\`enes des $p_j$ sont nulles sauf la $j$-i\`eme,
autrement dit que les $\cal{F}_j$ sont d\'efinis par les syst\`emes 
$$
n_{1,j}=\cdots = n_{r,j} = 0 \;\; (j=1,\ldots,n) , \;\;\; \sum_{\a=1}^n n_{1,\a}= \cdots = \sum_{\a=1}^n n_{r,\a}= 0.
$$
Des deux pr\'esentations obtenues de la m\^eme famille de feuilletages, on d\'eduit que 
$n_{a,j}=\sum_{b=1}^n C_{abj}m_{b,j}\,$ pour $j=1,\ldots,n$ et que 
$$
\sum_{b=1}^r\sum_{j=1}^n C_{abj}m_{b,j} = \sum_{b=1}^r\sum_{j=1}^n D_{ab}m_{b,j}.
$$
Donc $C_{abj}=D_{ab}$ pour tout $j$ et les deux bases 
de $1$-formes d\'efinissent la m\^eme structure presque grassmannienne.
\epf
Deux ${\rm G}_{r,n}$-structures, d\'efinies sur des 
germes $(M,x_0)$ et $(M',x'_0)$, sont isomorphes s'il existe 
un isomorphisme de $(M,x_0)$ sur $(M',x'_0)$ qui les \'echange.

Suivant la terminologie habituelle de la th\'eorie des $\rm G$-structures,
une ${\rm G}_{r,n}$-structure est {\em int\'egrable} 
si elle est (localement) d\'efinissable par une base de $1$-formes exactes.
Sans le dire, nous avons montr\'e au d\'ebut de la Section 1.2 que la grassmannienne
$\G_{r,n}$ admet une ${\rm G}_{r,n}$-structure int\'egrable.
On v\'erifie facilement que deux ${\rm G}_{r,n}$-structures int\'egrables 
sont isomorphes.

\sk
Il est remarquable que si $r\geq 2$, ce n'est pas le cas si $r=1$,
la propri\'et\'e pour un $d$-tissu de type $(r,n)$ d'\^etre isomorphe
\`a un tissu grassmannien est caract\'eris\'ee par 
une propri\'et\'e g\'eom\'etrique simple. Comme le remarque 
Akivis \cite{Ak}, si $d\geq n+1$ il suffit pour montrer que c'est le cas, 
de montrer d'une part que le tissu est presque grassmannien,
c'est-\`a-dire que tous les $(n+1)$-tissus qu'on peut en extraire
d\'efinissent la m\^eme structure presque grassmannienne,
ensuite que cette structure est int\'egrable.

Nous n'aurons pas directement besoin dans cet article 
des crit\`eres d'int\'egrabilit\'e donn\'es en particulier 
dans Akivis \cite{Ak} et Goldberg \cite{Go1}, mais 
ils interviennent dans \cite{P-Tr} dont nous  
utiliserons  les r\'esultats dans les d\'emonstrations 
du  Chapitre 4.

\ssct{Tissu tangent et applications}

Soit $\cal{T}=\{\cal{F}_1,\ldots,\cal{F}_d\}$ un $d$-tissu
de type $(r,n)$ sur un germe lisse $(M,x_0)$.
Pour  chaque $j$, soit $\O_j$ une normale g\'en\'eratrice de $\cal{F}_j$.

\sk
La valeur $\O_j(x_0)\in \wedge^r T^\star_{x_0}M$ est une normale 
g\'en\'eratrice d'un feuilletage constant $\cal{F}_{j,x_0}$ de $T_{x_0}M$,
\'evidemment ind\'ependant du choix de $\O_j$.
\bd
Le $d$-tissu constant $\cal{T}_{x_0} = \{\cal{F}_{1,x_0},\ldots,\cal{F}_{d,x_0}\}$
sur $T_{x_0}M$ est le tissu tangent de $\cal{T}$ en $x_0$.
\ed
Soit $z$ une fonction ou une forme sur $(M,x_0)$.
Sa {\em valuation} en $x_0$ est son ordre fini ou infini 
d'annulation au point $x_0$. Si la valuation de $z$ est sup\'erieure
ou \'egale \`a un entier donn\'e $h\geq 0$,
on lui associe sa {\em partie principale d'ordre $h$},
un polyn\^ome  $[z]_h$ homog\`ene de degr\'e $h$ sur $T_{x_0}M$.

Si l'on se donne un syst\`eme de coordonn\'ees de $M$
en $x_0$, $[z]_h$ est simplement la partie  homog\`ene de degr\'e $h$
de la s\'erie de Taylor de $z$ en $x_0$. Il est bien connu 
que, interpr\'et\'ee comme un polyn\^ome sur $T_{x_0}M$,
elle ne d\'epend pas du syst\`eme de coordonn\'ees choisi.

\bk
Soit $R$ l'espace des relations ab\'eliennes du tissu $\cal{T}$ en $x_0$
et si $h\in \N$, $R_h$ le sous-espace de $R$ de ses relations de valuation $\geq h$.

\sk
Si des $z_j\,\O_j$ sont les composantes d'un \'el\'ement de $R_h$,
comme 
$$
[z_j\,\O_j]_h = [z_j]_h\O_j(x_0), \;\;\; [d(z_j\,\O_j)]_{h-1} = d[z_j]_h \wedge \O_j(x_0),
$$
le $d$-uplet des $[z_j]_h\,\O_j(x_0)$ est une relation ab\'elienne, homog\`ene de degr\'e
$h$ du tissu tangent $\cal{T}_{x_0}$. Nous avons donc un morphisme lin\'eaire
\beq
\label{inj}
R_h/R_{h+1} \rightarrow R_{x_0}(h),
\eeq
\'evidemment  injectif, \`a valeurs dans l'espace $R_{x_0}(h)$ des relations ab\'eliennes homog\`enes 
de degr\'e $h$ du tissu tangent $\cal{T}_{x_0}$. Compte tenu de la Proposition \ref{CG-cste} nous obtenons la borne de Chern-Griffiths 
sous la forme pr\'ecis\'ee
\beq
\label{borne-precise}
\dim R_h/R_{h+1} \leq   \dim E_r(h)\times  \max\,(d-(r+h)(n-1)-1,0).
\eeq

\sk
Nous supposons maintenant que le tissu $\cal{T}$ est semi-extr\'emal
de type $(r,n)$ avec $r\geq 2$.

\sk
Les applications canoniques $R_0/R_1\rightarrow R_{x_0}(0)$
et $R_1/R_2\rightarrow R_{x_0}(1)$ sont alors n\'ecessairement des isomorphismes,
en particulier le tissu tangent $\cal{T}_{x_0}$ 
est semi-extr\'emal. Comme la propri\'et\'e d'\^etre semi-extr\'emal
est ouverte, le tissu tangent $\cal{T}_x$ en tout point $x\in (M,x_0)$
est semi-extr\'emal. 

\sk
Nous pouvons donc appliquer la Proposition \ref{tensor}. La 
notion de structure presque grassmannienne permet d'\'enoncer le r\'esultat 
dans un cadre g\'eom\'etrique naturel. 
\bt
\label{forme-norm}
Un $d$-tissu semi-extr\'emal $\cal{T}= \{\cal{F}_1,\ldots,\cal{F}_d\}$ de type $(r,n)$, $r\geq 2$,
sur un germe $(M,x_0)$, y d\'etermine une ${\rm G}_{r,n}$-structure 
pour laquelle il est presque grassmannien. De plus,
si la base $(m_{a,\a})_{(a,\a)\in A}$ d\'efinit la structure,
il existe une base $(p_\a(t))_{\a=1}^n$ du module des polyn\^omes
de degr\'e $\leq n-1$ et \`a coefficients analytiques sur $(M,x_0)$, et des fonctions 
$\t_j: (M,x_0)\rightarrow  \P^1$ 
telles que $\cal{F}_j$ est d\'efini
par le syst\`eme 
$$
\sum_{\a=1}^n p_\a(\t_j) \, m_{a,\a}=0, \;\; a =1,\ldots,r.
$$
\et
Il r\'esulte  des commentaires de la section pr\'ec\'edente 
que le tissu est isomorphe \`a un germe de tissu grassmannien si et seulement
si l'on peut {\em choisir} les formes $m_{a,\a}$  ferm\'ees.

\sk
Un changement de base $n_{a,\a} = \sum_{\b=1}^n C_{\a\b}m_{a,\b}$
permet de supposer $p_\a(t)=t^{\a-1}$, ind\'ependant de $x\in (M,x_0)$. 
On se ram\`ene aussi \`a ce que les $\t_j$ soient \`a valeurs 
dans $\C$. On obtient ainsi {\em la forme normale} suivante 
\beq
\label{norm1}
\sum_{\a=1}^n\t_j^{\a-1} \, m_{a,\a}=0, \;\;\;\;  a =1,\ldots,r,
\eeq
des \'equations des $\cal{F}_j$. C'est cette forme que nous utiliserons 
dans les calculs qui vont suivre. Les points $[1:\t_j(x):\cdots:\t_j(x)^{n-1}]$
appartiennent \`a une courbe rationnelle normale de degr\'e $n-1$ de $\P^{n-1}$.
C'est une propri\'et\'e des normales du tissu qui, sauf si $n=2$,
n'est pas impliqu\'ee par le fait qu'il est presque grassmannien.
Enfin, m\^eme si la structure est int\'egrable, on ne peut pas 
toujours choisir une base de formes ferm\'ees pour laquelle 
les \'equations du tissu s'\'ecrivent sous la forme normale

\ssct{Une d\'emonstration historique}

Soit $\cal{T}$ un $d$-tissu de type $(r,n)$ avec $q(d)=n-1$.
Compte tenu de (\ref{borne-precise}) il n'a pas
de relation ab\'elienne non triviale de valuation $\geq 2$.
Il revient au m\^eme de dire qu'il est de rang maximal 
ou qu'il est semi-extr\'emal, ce que nous supposons.

\sk
C'est un cas qui n'est pas exclu dans la suite de cet article,
mais il serait dommage de ne pas rappeler ici la d\'emonstration
classique du fait qu'il est isomorphe \`a un germe 
de tissu alg\'ebrique grassmannien.
Nous renvoyons \`a Chern \cite{C} pour sa relation avec le th\'eor\`eme de Lie 
sur les surfaces de double translation et 
ses g\'en\'eralisations (Tschebotarow 1927, Wirtinger 1938).

\sk
La puissance de la m\'ethode canonique
est spectaculaire dans ce cas. La d\'emonstration de Howe lorsque $(r,n)=(1,2)$ 
(voir \cite{BB}) se g\'en\'eralise sans encombre, \`a partir \og seulement \fg\,
des d\'efinitions, de la borne de Chern-Griffiths et 
du th\'eor\`eme d'Abel inverse. 

\sk
Le $d$-tissu $\cal{T}$ est de rang maximal $r+n$. 
Nous reprenons les notations de la Section 1.3.
Les applications de Poincar\'e prennent leurs valeurs dans 
l'espace $\P^{r+n-1}$. C'est un exercice de v\'erifier directement
(plus simplement que dans la prochaine section) que si $x\in (M,x_0)$, 
les $d$ points $\k_j(x)$ engendrent un espace projectif $\lan \k_1(x), \ldots, \k_d(x) \ran$
de dimension $n-1$ et que l'application 
$$
x\mapsto \lan \k_1(x), \ldots, \k_d(x) \ran,
$$
est un isomorphisme de $(M,x_0)$ sur un germe de la grassmannienne $\G_{r,n}$.
C'est en fait cette application que Poincar\'e consid\`ere dans \cite{Po}
avec $r=1$ et $n=2$.

On v\'erifie aussi facilement que, par construction,  l'image du tissu $\cal{T}$ par 
cet isomorphisme est un tissu grassmannien. Cette image 
est un tissu alg\'ebrique compte tenu du th\'eor\`eme
d'Abel inverse.

\ssct{Mise en oeuvre de la m\'ethode canonique}

Nous continuons l'\'etude d'un $d$-tissu semi-extr\'emal 
$\cal{T}= \{\cal{F}_1,\ldots,\cal{F}_d\}$, de type $(r,n)$ avec $r\geq 2$,
sur un germe lisse $(M,x_0)$. 

Les propri\'et\'es d\'ej\`a obtenues 
sont r\'esum\'ees par le Th\'eor\`eme \ref{forme-norm}.
Nous choisissons la base de $1$-formes $(m_{a,\a})_{(a,\a)\in A}$ 
telle que les feuilletages $\cal{F}_j$ sont pr\'esent\'es 
sous la forme normale (\ref{norm1}).
La $r$-forme  
\beq
\label{normales}
\O_j
 = \bigwedge_{a=1}^r \left(\sum_{\a=1}^n \t_j^{\a-1} m_{a,\a}\right) =  \sum_{\rho=0}^{r(n-1)} \t_j^\rho \, K_\rho,
\eeq
o\`u le membre de droite est simplement le d\'eveloppement du membre 
m\'ediant en puissance de $\t_j$, est une normale g\'en\'eratrice
du feuilletage $\cal{F}_j$. Elle n'est pas ferm\'ee en g\'en\'eral. 
Les $r$-formes $K_\rho$ sur $(M,x_0)$
sont lin\'eairement ind\'ependantes au point $x_0$, puisque 
le syst\`eme des normales en un point est de rang $r(n-1)+1$,

\sk
Un $d$-uplet de normales $z_j\,\O_j$, {\em suppos\'ees ferm\'ees}, est donc une 
relation ab\'elienne si et seulement si les fonctions $z_j$ v\'erifient les \'equations 
\beq
\label{V0}
\sum_{j=1}^d \t_j^\rho \, z_j =0, \;\;\;\; \rho=0,\ldots,r(n-1).
\eeq
Le syst\`eme~(\ref{V0}) est un syst\`eme de Vandermonde sous-d\'etermin\'e 
dont la solution g\'en\'erale est facile
\`a \'ecrire \`a partir de la solution $(c_1,\ldots,c_d)$
du syst\`eme de Vandermonde-Cramer
\beq
\label{V1}
\sum_{j=1}^d \t_j^\rho c_j = \delta_{\rho \,d-1}, \;\;\; \rho=0,\ldots,d-1.
\eeq
\ble
\label{V}
La solution (analytique) g\'en\'erale du syst\`eme (\ref{V0}) sur $(M,x_0)$
est donn\'ee par 
\beq
\label{V2}
z_j(x)=c_j(x)f(x,\t_j(x)), \;\;\;\; j=1,\ldots,d,
\eeq
o\`u $f(x,t)$ est (analytique et) polynomiale en $t$ de degr\'e $\leq q(d)$.
\ele
Rappelons que nous notons $q(d)=d-r(n-1)-2$. 
Insistons d'autre part sur le fait que (\ref{V2}) caract\'erise les syst\`emes de normales 
$z_j\O_j$ {\em non n\'ecessairement ferm\'ees}, de somme nulle.

\bk
Soit $N+1$ le rang du tissu $\cal{T}$. 
Suivant la m\'ethode canonique, nous choisissons une base $B$,
$$
(\phi_1^{(\l)}, \ldots , \phi_d^{(\l)}), \;\;\;  \l=1,\ldots,N+1,
$$
de l'espace $R$ des relations ab\'eliennes du tissu $\cal{T}$
et nous introduisons les applications de Poincar\'e, une pour chaque feuilletage $\cal{F}_j$,
\beq
\label{Po}
{\k}_j : \; (M,x_0) \rightarrow \P^N, \qquad {\k}_j (x) = [\phi_j^{(1)}(x): \cdots : \phi_j^{(N+1)}(x)].
\eeq
Nous verrons  bient\^ot qu'au moins une composante $\phi_j^{(\l)}(x)$ est non nulle.

\sk
L'\'ecriture (\ref{Po}) est un peu abusive : les coordonn\'ees homog\`enes de ${\k}_j(x)$
sont des formes, mais colin\'eaires. Elle a l'int\'er\^et de  
montrer que ces applications ne supposent pas un choix des normales 
g\'en\'eratrices et que, \`a un automorphisme
pr\`es de l'espace but $\P^N$, la famille des applications ${\k}_j$
ne d\'epend pas de la base $B$ choisie.

\sk
Si l'on \'ecrit $\phi_j^{(\l)}=\zeta_j^{(\l)}N_j$ o\`u les $N_j$ sont des normales 
g\'en\'eratrices {\em ferm\'ees}, on obtient $d\zeta^{(\l)}_j\wedge N_j=0$. Autrement dit
$\k_j=[\zeta_j^{(1)}: \cdots : \zeta_j^{(N+1)}]$ est 
constante le long des feuilles de  $\cal{F}_j$.

\sk
Remarquons aussi que, lorsque $x$ varie dans $(M,x_0)$, 
les points $\k_j(x)$ engendrent $\P^N$. En effet, s'ils 
restaient contenus dans un hyperplan, on aurait une relation 
de d\'ependance lin\'eaire non triviale entre les relations 
ab\'eliennes $(\phi_1^{(\l)},\ldots,\phi^{(\l)}_d)$, en contradiction
avec l'hypoth\`ese qu'elles forment une base.

\bk
Nous utilisons les normales (\ref{normales})
et nous \'ecrivons $\phi_j^{(\l)} = z_j^{(\l)} \, \O_j$. Le point 
${\k}_j (x)$ s'\'ecrit 
$$
\k_j(x)= [z_j^{(1)}(x): \cdots : z_j^{(N+1)}(x)].
$$
Les propri\'et\'es que nous voulons \'etablir, d'abord au point $x_0$,
ne d\'ependent pas du choix de la base $B$. Nous allons la choisir de telle sorte que 
ces propri\'et\'es apparaissent comme \'evidentes.

\sk
Nous pouvons supposer $(M,x_0)=(\C^{rn},0)$. Notons $u_{a,\a}$ les
formes lin\'eaires sur $\C^{rn}$ telles que $du_{a,\a} = m_{a,\a}(0)$.
Les  feuilletages tangents $\cal{F}_{j,x_0}$ des $\cal{F}_j$
sont d\'efinis 
par  les syst\`emes $dl_{j,1}=\cdots=dl_{j,r}=0$, o\`u  
$$
l_{j,a}(x) =  \sum_{\a=1}^n \t_j(x_0)^{\a-1} \, u_{a,\a}(x).
$$

\sk
Nous choisissons la base $B$ de $R$ en partant d'un syst\`eme de $q(d)+1$
relations de valuation $0$ dont les parties principales constituent la base
$$
(c_1(x_0)\t_1(x_0)^\rho \,\O_1(x_0), \ldots ,c_d(x_0)\t_d(x_0)^\rho \, \O_d(x_0)), \;\;\;   \rho = 0,\ldots,q(d),
$$
de l'espace des relations homog\`enes de degr\'e $0$ de $\cal{T}_{x_0}$
et en la compl\'etant par une base $B_1$ de l'espace 
$R_1$ des relations de valuation $\geq 1$ de $\cal{T}$.

\sk
Les $c_j(x_0)$ sont non nuls compte tenu de la forme du syst\`eme~(\ref{V1})
donc 
$$
{\k}_j(x_0) = [1 : \t_j(x_0) : \cdots : \t_j(x_0)^{q(d)} : 0 : \cdots : 0],
$$ 
avec $\t_j(x_0)\neq \t_k(x_0)$ si $j\neq k$.
Les applications $\k_j$ sont bien d\'efinies,
les points $\k_j(x_0)$ engendrent un $q(d)$-plan $\,\P^{q(d)}(x_0)$
et sont en position g\'en\'erale dans cet espace.
Ils appartiennent \`a la courbe rationnelle normale $\s(x_0)$ de degr\'e 
$q(d)$ param\'etr\'ee par $t\mapsto [1 : t : \cdots : t^{q(d)} : 0 : \cdots :0]$.
C'est la seule courbe de ce type (une telle courbe est d\'etermin\'ee
par $q(d)+3$ de ses points) qui passe par les points $\k_j(x_0)$.

\bk
Nous choisissons la base $B_1$ en compl\'etant un syst\`eme 
de $r$ relations de valuation $1$ dont les parties principales 
sont les relations ab\'eliennes homog\`enes de degr\'e $1$
suivantes du tissu tangent $\cal{T}_{x_0}$ : 
$$
(c_1(x_0) \,l_{1,a}(x) \O_1(x_0), \ldots ,c_d(x_0)l_{d,a}(x)\O_d(x_0)), \;\;\; 1\leq a\leq r.
$$
Ce choix fait, parmi les composantes de valuation $\geq 1$ de $\k_j(x)$
figurent 
$$
l_{j,1}(x)+ O(|x|^2), \ldots, l_{j,r}(x)+ O(|x|^2),
$$
dont les diff\'erentielles en $x_0$ sont lin\'eairement ind\'ependantes.
Il en r\'esulte que ${\k}_j$ est de rang $\geq r$ en $x_0$,
donc de rang $r$ puisque $\k_j$ est constante le long 
des feuilles de $\cal{F}_j$. Son image est un germe $(Z,\k_j(x_0))$
de vari\'et\'e lisse de dimension $r$, transverse au plan $\P^{q(d)}$ 
et donc \`a la courbe $\s(x_0)$.

Enfin, comme les composantes homog\`enes d'un point de $\P^{q(d)}(x_0)$, correspondant aux composantes 
pr\'ec\'edentes de $\k_j(x)$, sont nulles, les formes $dl_{j,a}(x_0)$
s'annulent sur le noyau de la diff\'erentielle de l'application 
$x\mapsto \P^{q(d)}(x_0)$ pour tout $j$ et tout $a$.
Comme les feuilletages $\cal{F}_j$ v\'erifie la condition (PG), ce noyau est nul. 

Ainsi $x\mapsto \P^{q(d)}(x)$ est une immersion et comme 
$\P^{q(d)}(x)$ est l'espace engendr\'e par la courbe $\s(x)$, 
l'application $x\mapsto \s(x)$ en est une aussi.

\sk
Les propri\'et\'es pr\'ec\'edentes sont v\'erifi\'ees en tout point : 
il suffit d'appliquer la m\^eme r\'eduction
\`a la restriction pr\`es de $x\in (M,x_0)$ 
de l'espace des relations ab\'eliennes du tissu $\cal{T}$
\footnote{En fait le rang de $\cal{T}$ en un point est constant
mais c'est sans importance ici.}.
Nous avons donc :
\ble
\label{appli-P}
Chaque application de Poincar\'e ${\k}_j: (M,x_0) \rightarrow \P^N$ 
est de rang $r$, constante le long des feuilles de $\cal{F}_j$.
La r\'eunion des germes d'image des $\k_j$ engendre $\P^N$. 

Pour tout $x\in (M,x_0)$,
les $d$ points ${\k}_j(x)$ sont des points deux-\`a-deux distincts 
d'une courbe rationnelle normale $\s(x)$ de degr\'e $q(d)$,
transverse aux images des germes $\k_j$.
L'application $x\mapsto\s(x)$ est une immersion.
\ele

\ssct{Conditions d'int\'egrabilit\'e}

Cette section est une section de calculs. Il s'agit,
\'etant donn\'e {\em un seul} feuilletage d\'efini par un syst\`eme 
{\em compl\`etement int\'egrable} de la forme (\ref{norm1})
avec une fonction $\t$ donn\'ee,
d'\'etablir certaines formules en gardant trace de leur 
d\'ependance en $\t$. Elles seront utiles 
dans la prochaine section.

\sk
Soit $(X_{a,\a})_{(a,\a)\in A}$ la base de champs de vecteurs duale 
de la base de $1$-formes $(m_{a,\a})_{(a,\a)\in A}$. 
Pour $t\in \C$, nous posons  
$$
l_a(t) = \sum_{\a=1}^n t^{\a-1}m_{a,\a}, \;\; a=1,\ldots,r, \qquad \O(t) = \bigwedge_{a=1}^r l_a(t).
$$
Pour $t$ fix\'e, les formes $l_a(t)$ et les formes $m_{a,\a}$ avec $\a\neq 1$ constituent 
une base de $1$-formes sur $(M,x_0)$ dont la base duale est constitu\'ee des champs de vecteurs 
$X_{a,1}$ (ind\'ependants de $t$)  et
$$
Y_{a,\a}(t) = X_{a,\a} - t^{\a-1}X_{a,1}, \;\; a=1,\ldots,r, \;\; \a=2,\ldots,n.
$$
Dans la suite comme ci-dessus, nous ne notons pas la variable $x$.
En particulier, si $\t$ est une fonction sur $(M,x_0)$, nous notons 
$l_a(\t)$ pour la $1$-forme $x\mapsto l_a(x,\t(x))$ et $Y_{a,\a}(\t)$
pour le champ $x\mapsto Y_{a,\a}(x,\t(x))$.
\ble
\label{calcul1}
Il existe, pour tout $(a,\a)\in A$ avec $\a\neq 1$, des polyn\^omes $A_{a,\a}(t)$ de degr\'es $\leq 2n-1$
et $B_{a,\a}(t)$ de degr\'es $\leq 2n-2$, \`a coefficients analytiques 
sur $(M,x_0)$, tels qu'on ait les propri\'et\'es suivantes.

Pour toute fonction $\t: (M,x_0)\rightarrow \C$ d\'efinissant 
un syst\`eme diff\'erentiel $l_a(\t)=0$, $\,a=1,\ldots,r$,
compl\`etement int\'egrable, et pour toute solution $z$ sur $(M,x_0)$
de l'\'equation $d(z\O(\t))=0$, on a pour $(a,\a)\in A$ avec $\a\neq 1$ : 
\beq
\label{Y-theta}
Y_{a,\a}(\t)\cdot \t = A_{a,\a}(\t),
\eeq
\beq
\label{Y-z}
Y_{a,\a}(\t)\cdot z   =  zX_{a, 1}\cdot \t^{\a-1} + zB_{a,\a}(\t).
\eeq
\ele
\bpf
Il est commode de noter $[k]_I(t)$, o\`u $k\in \N$ et o\`u $I$ est un indice qui d\'epend du contexte,
une forme ou une fonction sur $(M,x_0)$, polynomiale de degr\'e $\leq k$
en $t\in \C$. Par exemple on peut \'ecrire $dm_{a,\a} = [0]_{a,\a}(t)$ et,
de la d\'efinition $l_a(\t)=\sum_{\a=1}^n \t^{\a-1} m_{a,\a}$, on d\'eduit   
$$
dl_a(\t) 
 = 
\sum_{\b=1}^n   d\t^{\b-1} \wedge m_{a,\b} + [n-1]_{a,\a}(\t).
$$
D'autre part, le syst\`eme \'etant suppos\'e compl\`etement int\'egrable, il existe des $1$-formes $\Gamma_{a,b}^\t$,
qui d\'ependent de $\t$ d'une fa\c{c}on inconnue, telles que 
$$
dl_a(\t) 
= \sum_{b=1}^r  \Gamma_{a,b}^\t \wedge l_b(\t)
= \sum_{b=1}^r \sum_{\b=1}^n \t^{\b-1} \Gamma_{a,b}^\t \wedge m_{b,\b}.
$$
La forme  $dl_a(\t)$ a aussi une d\'ecomposition  canonique de la forme 
$\sum C_{b,\b,c,\g} m_{b,\b} \wedge m_{c,\g}$ avec $C_{b,\b,c,\g} + C_{c,\g,b,\b} = 0$.
Nous chassons le coefficient de $m_{a,\a}\wedge m_{b,\b}$ dans les deux expressions 
ci-dessus de $dl_a(\t)$.

\bk
\nk
1) En identifiant les coefficients de $m_{a,\a}\wedge m_{b,\b}$ {avec $a\neq b$},
on obtient :  
$$
\lan d\t^{\a-1}, X_{b,\b}\ran  + [n-1]_{a,b,\b}(\t) 
= 
\t^{\a-1} \lan\Gamma_{a,a}^\t,X_{b,\b}\ran  - \t^{\b-1}\lan\Gamma_{a,b}^\t,X_{a,\a}\ran,
$$
soit si $\b=1$ :
$$
\lan d\t^{\a-1}, X_{b,1}\ran  + [n-1]_{a,b,1}(\t) 
= 
\t^{\a-1} \lan\Gamma_{a,a}^\t,X_{b,1}\ran  - \lan\Gamma_{a,b}^\t,X_{a,\a}\ran.
$$
\'Eliminons $\Gamma_{a,b}^\t$ entre les deux \'equations pr\'ec\'edentes :
$$
\lan d\t^{\a-1}, Y_{b,\b}(\t)\ran  + [2n-2]_{a,b,\b}  =  \t^{\a-1} \lan\Gamma_{a,a}^\t,Y_{b,\b}(\t)\ran.
$$
En choisissant $\a=1$ puis $\a=2$, nous obtenons pour $a\neq b$ :
$$
\lan\Gamma_{a,a}^\t,Y_{b,\b}(\t)\ran = [2n-2]_{a,b,\b}(\t),
$$
$$
\lan d\t, Y_{b,\b}(\t)\ran  =  \t\lan\Gamma_{a,a}^\t,Y_{b,\b}(\t)\ran + [2n-2]_{a,b,\b}(\t).
$$

\sk
{\em Parce que  $r\geq 2$, nous pouvons effectivement choisir $a\neq b$
si $b$ est donn\'e.}
La seconde relation, compte tenu de la premi\`ere, 
donne alors $\lan d\t, Y_{b,\b}(\t)\ran = [2n-1]_{b,\b}(\t)$. C'est la formule (\ref{Y-theta}).

\sk
D'autre part, la premi\`ere relation, valable si $a\neq b$, donne :
\beq
\label{A}
\lan \sum_{a=1}^r \Gamma_{a,a}^\t,Y_{b,\b}(\t)\ran = \lan \Gamma_{b,b}^\t,Y_{b,\b}(\t)\ran + [2n-2]_{b,\b}(\t).
\eeq

\bk
\nk
2) En identifiant les coefficients de $m_{a,\a}\wedge m_{a,1}$ pour $\a\neq 1$,
on obtient : 
$$
\lan d\t^{\a-1},X_{a,1}\ran +  [n-1]_{a,\a}(\t) =
\t^{\a-1} \lan \Gamma_{a,a}^\t,X_{a,1}\ran 
 - 
\lan \Gamma_{a,a}^\t,X_{a,\a}\ran \, ;
$$
c'est-\`a-dire 
$$
 \lan \Gamma_{a,a}^\t,Y_{a,\a}(\t)\ran  = - \lan d\t^{\a-1},X_{a,1}\ran +  [n-1]_{a,\a}(\t)).
$$
Compte tenu de (\ref{A}), on a donc :
\beq
\label{C}
\lan \sum_{a=1}^r \Gamma_{a,a}^\t,Y_{b,\b}(\t)\ran = - \lan d\t^{\b-1},X_{b,1}\ran  + [2n-2]_{b,\b}(\t).
\eeq

\bk
\nk
3) La d\'eriv\'ee ext\'erieure $d\O(\t)$ de $\O(\t) = \wedge_{a=1}^r l_a(\t)$ est 
\'egale \`a 
$$
\sum_{b=1}^r (-1)^{b-1} l_1\t)\wedge \cdots \wedge l_{b-1}(\t)\wedge dl_b(\t) \wedge l_{b+1}(\t)\wedge \cdots \wedge l_r(\t)
$$
donc \`a $(\sum_{b=1}^r \Gamma^\t_{b,b}) \wedge \O(\t)$.

Soit maintenant $z$ une solution sur  $(M,x_0)$ de l'\'equation $d(z\O(\t))=0$ :
$$
(dz + z\sum_{b=1}^r \Gamma^\t_{b,b}) \wedge \O(\t)=0,
$$
donc $\lan Y_{a,\a}(\t), dz + z\sum_{b=1}^r \Gamma^\t_{b,b} \ran =0$ pour tout $\a\neq 1$.
Compte tenu de la formule (\ref{C}), nous obtenons 
$$
Y_{a,\a}(\t)\cdot z = zX_{a,1}\cdot \t^{\a-1}  + z[2n-2]_{a,\a}(\t), 
$$
ce qui termine la d\'emonstration du lemme.
\epf
\ble
\label{calcul2}
Sous les hypoth\`eses du lemme pr\'ec\'edent, il existe des polyn\^omes 
$a_{a,\a}(t),\, b_{a,\a}(t)$ en $t$ ainsi que des polyn\^omes $c_{a,\a}(t,s)$
de degr\'e $\leq 2n-3$ et $d_{a,\a}(t,s)$ de degr\'e $\leq n-2$ en $(t,s)$,
\`a coefficients analytiques et ind\'ependants de la fonction $\t$, tels
que :
$$
Y_{a,\a}(t)\cdot \left(\fr{z}{t-\t}\right) 
=
z\left(\fr{a_{a,\a}(t)}{t-\t} + \fr{b_{a,\a}(t)}{(t-\t)^2} + c_{a,\a}(t,\t)\right)  + X_{a1}\cdot (zd_{a,\a}(t,\t)),
$$
\ele
\bpf
Nous ne notons pas $a$ et $\a$ qui sont fix\'es. Nous posons  
$u = z/(t-\t)$, $\,v = z/(t-\t)^2$. Compte tenu de (\ref{Y-theta})--(\ref{Y-z}), il vient : 
$$
Y(\t)\cdot u = \fr{Y(\t)\cdot z}{t-\t} + z\, \fr{Y(\t)\cdot \t}{(t-\t)^2}  
= 
u(X_1\cdot \t^{\a-1} + B(\t))  + vA(\t).  
$$
En utilisant les d\'eveloppements 
$$
A(s) = A(t) - (t-s)A'(t) + (t-s)^2 \tilde{A}(t,s), \;\; B(s) = B(t) - (t-s)\tilde{B}(t,s),
$$
o\`u $\tilde{A}(t,s)$ et $\tilde{B}(t,s)$ sont de degr\'e $\leq 2n-3$ en $(t,s)$,
nous obtenons  
$$
Y(\t)\cdot u = u(B(t)-A'(t)) + vA(t) + uX_1\cdot \t^{\a-1} + zC(t,\t),
$$
o\`u $C(t,s)$ est polynomial en $(t,s)$ de degr\'e $\leq 2n-3$.

\sk
D'autre part, $\,Y(t)\cdot u = Y(\t)\cdot u  - (t^{\a-1}-\t^{\a-1})X_1\cdot u$,
et en \'ecrivant  $t^{\a-1}-s^{\a-1} = (t-s)S(t,s)$, nous obtenons  
$$
Y(t)\cdot u = Y(\t)\cdot u  - X_1\cdot ((t-s)S(t,\t)u) - uX_1\cdot \t^{\a-1},
$$
o\`u $S(t,s)$ est de degr\'e $\leq n-2$ en $(t,s)$. Quand on remplace 
$Y(\t)\cdot u$ par l'expression que nous en avons obtenue ci-dessus, les
termes en  $X_1\cdot \t^{\a-1}$ s'\'eliminent, ce qui donne le lemme.
\epf

\ssct{L'application de Blaschke}

Nous poursuivons maintenant la discussion amorc\'ee dans la Section 3.4, dont 
nous reprenons les notations.
Le tissu $\cal{T}$ est pr\'esent\'e sous la forme normale (\ref{norm1}).
Ses normales g\'en\'eratrices $\O_j$ sont d\'efinies par (\ref{normales}).
Nous nous donnons une  base $B$ de l'espace de ses relations ab\'eliennes.
Les applications de Poincar\'e $\k_j$ sont d\'efinies par (\ref{Po}). Leurs propri\'et\'es 
sont d\'ecrites dans le Lemme \ref{appli-P}. En particulier, pour 
$x$ donn\'e, les points $\k_j(x)$ appartiennent \`a une 
courbe rationnelle normale $\s(x)$, de degr\'e $q(d)$.

\sk
Posons
$$
P(x,t) = \prod_{j=1}^d (t-\t_j(x)).
$$
Il s'agit d'abord de param\'etrer la courbe $\s(x)$ qui passe par 
les points $\k_j(x) = [\phi_j^{(1)}(x): \cdots : \phi_j^{(N+1)}(x)]$.
\'Ecrivons $\phi_j^{(\l)}=z_j^{(\l)}\,\O_j$ et 
$$
z_j(x)=(z_j^{(1)}(x),\ldots,z_j^{(N+1)}(x))\in \C^{N+1}\bck\{0\}.
$$
\ble
Pour tout $x\in (M,x_0)$, la courbe $\s(x)$ est la projection canonique dans $\P^N$
de la courbe 
de $\C^{N+1}\bck\{0\}$ param\'etr\'ee par 
\beq
\label{z}
z(x,t)  = \sum_{j=1}^d  \fr{P(x,t)}{t-\t_j(x)} \, z_j(x).
\eeq
\ele
\bpf
Le point $x$ de $(M,x_0)$ ne joue aucun r\^ole dans ce qui suit, nous ne le notons pas.

L'application $t\mapsto z(t)$ est polynomiale de degr\'e $\leq d-1$ et la projection 
canonique de son image passe par les points $\k_j$.
Il suffit de montrer que $z(t)$ est en fait de degr\'e $\leq q(d)$.

\sk
C'est une cons\'equence du Lemme \ref{V}. Les composantes de $z_j(t)$
sont de la forme $z_j^{(\l)}(t) = c_jf^{(\l)}(\t_j)$ o\`u $f^{(\l)}(t)$ est 
un polyn\^ome de degr\'e $\leq q(d)$ et $(c_1,\ldots,c_d)$
la solution du syst\`eme (\ref{V1}). 
Montrons plus g\'en\'eralement que si $f(t)$ est un polyn\^ome de degr\'e $\leq d-1$
alors 
$$
\sum_{j=1}^d  \fr{P(t)}{t-\t_j} \, c_jf(\t_j) = f(t),
$$
ce qui donne le r\'esultat. 
Compte tenu de (\ref{V1}), on a d'abord 
$$
\sum_{j=1}^d  \fr{P(t)}{t-\t_j} \, c_j = \sum_{j=1}^d  \fr{P(t)-P(\t_j)}{t-\t_j} c_j= 1,
$$
puisque $P(\t_j)=0$ et que pour $t$ fix\'e, $(P(t)-P(s))/(t-s)$
est un polyn\^ome unitaire de degr\'e $d-1$ en $s$.
En g\'en\'eral on en d\'eduit 
$$
\sum_{j=1}^d  \fr{P(t)}{t-\t_j} \, c_jf(\t_j) = f(t) - P(t) \sum_{j=1}^d  \fr{f(t)-f(\t_j)}{t-\t_j} \, c_j = f(t),
$$
si $f(t)$ est de degr\'e $\leq d-1$ car pour $t$ fix\'e, 
$(f(t)-f(s))/(t-s)$ est un polyn\^ome de degr\'e $\leq d-2$ en $s$.
\epf
L'application 
\beq
\label{kappa}
\k: (M,x_0) \times \P^1 \rightarrow \P^N
\eeq
est d\'efinie en composant l'application $z\mapsto z(x,t)$
avec  la projection canonique. 
Elle est de rang $\geq r+1$ en tout point. 
En effet, en choisissant la  base $B$ comme dans la Section~3.4, dont on reprend 
les notations, on trouve que $t + O(|x|)$ et
$\sum_{\a=1}^n t^{\a-1} x_{a,\a} + O(|x|^2)$  pour $a=1,\ldots,r$, figurent dans un syst\`eme
de coordonn\'ees homog\`enes de $\k(x,t)$.

\bk
Le r\'esultat crucial suivant nous permettra dans le prochain chapitre
de d\'ecrire pr\'ecis\'ement la vari\'et\'e de Blaschke du tissu $\cal{T}$.
\bpr
\label{fond}
L'application (\ref{kappa}) est de rang constant $r+1$.
Le noyau de sa diff\'erentielle au point $(x,t)\in (M,x_0)\times \P^1$ est donn\'e par le syst\`eme 
\beq
\label{noyau}
dt + J(x,t)=0, \qquad \sum_{\a=1}^n t^{\a-1} m_{a,\a}(x)=0, \;\; a=1,\ldots,r, 
\eeq
o\`u $J(x,t)$ est une $1$-forme sur $(M,x_0)$ d\'ependant analytiquement de $t$.
\epr
\bpf
Comme nous savons que $\k$ est de rang $\geq r+1$, il reste \`a montrer que $d\k(x,t)$
s'annule sur l'espace d\'efini par (\ref{noyau}).

Si l'on exprime localement $\k$ dans une carte affine de $\P^N$,
cela revient \`a dire que $d\k$
est de la forme $d\k = C_0(dt+J) + \sum_{a=1}^r C_al_a(t)$.
Par d\'efinition des champs de vecteurs $Y_{a,\a}(t)$ pour $\a\neq 1$, nous avons :
$$
d\k = \fr{\pl \k}{\pl t} dt + \sum_{a=1}^r\sum_{\a=2}^n (Y_{a,\a}(t)\cdot \k)m_{a,\a}+ \sum_{a=1}^r (X_{a1}\cdot \k) l_a(t).
$$
Il s'agit donc de montrer, toujours dans une carte affine de $\P^N$,
qu'on a des formules de la forme 
\beq
\label{E}
Y_{a,\a}(t)\cdot \k = g_{a,\a}\fr{\pl \k}{\pl t}, \qquad \a\neq 1,
\eeq
o\`u les $g_{a,\a}$ sont des fonctions sur $(M,x_0)\times \C$. 

\sk
Il suffit de d\'emontrer ces formules  pour $t\notin \{\t_1(x),\ldots,\t_d(x)\}$.
Nous pouvons alors prendre les 
$$
z^{(\l)}(t)  = \sum_{j=1}^d  \fr{z_j^{(\l)}}{t-\t_j}, \qquad  \l=1,\ldots,N+1.
$$
comme syst\`eme de coordonn\'ees homog\`enes de $\k$. Les $z_j^{(\l)}$
sont des fonctions sur $(M,x_0)$ qui v\'erifient (\ref{V0}) 
et, pour tout $j$, $d(z_j^{(\l)}\O_j)=0$.

\sk
Comme les syst\`emes qui d\'efinissent les feuilletages $\cal{F}_j$ 
sont compl\`e\-tement int\'egrables, nous pouvons appliquer le Lemme \ref{calcul2}
\`a chacune des fonctions $z_j^{(\l)}$ ce qui nous donne 
$$
Y_{a,\a}(t)\cdot \fr{z_j^{(\l)}}{t-\t_j}
 =
z_j^{(\l)} \left(\fr{a_{a,\a}(t)}{t-\t_j} + \fr{ b_{a,\a}(t)}{(t-\t_j)^2} +  c_{a,\a}(t,\t_j) \right) 
+  X_1\cdot (z_j^{(\l)}d_{a,\a}(t,\t_j)). 
$$
En sommant par rapport \`a $j$, nous obtenons 
$$
Y_{a,\a}(t)\cdot z^{(\l)} = a_{a,\a}(t)z^{(\l)} - b_{a,\a}(t)\fr{\pl z^{(\l)}}{\pl t}
$$
car  $c_{a,\a}(t,s)$ et $d_{a,\a}(t,s)$ sont  de degr\'es $\leq 2n-1$ en $(t,s)$, donc  
$$
\sum_{j=1}^d z_j^{(\l)} c_{a,\a}(t,\t_j) +  X_1\cdot (\sum_{j=1}^d z_j^{(\l)}d_{a,\a}(t,\t_j)) = 0,
$$
compte tenu de (\ref{V0}). Au voisinage d'un point $(x,t)$, 
on choisit $\l_{0}$ tel que $z^{(\l_0)}\neq 0$. Les $z^{(\l)}/z^{(\l_0)}$ avec $\l\neq \l_0$
forment un syst\`eme de coordonn\'ees affines de $\k$ et 
$$
(Y_{a,\a}(t) - b_{a,\a}(t)\fr{\pl}{\pl t}) \cdot \fr{z^{(\l)}}{z^{(\l_0)}} = 0,
$$
ce qui termine la d\'emonstration.
\epf

\sct{R\'esultats d'alg\'ebricit\'e}

Dans cette section nous d\'emontrons nos principaux r\'esultats.
D'abord, nous montrons que la vari\'et\'e de Blaschke d'un $d$-tissu semi-extr\'emal
de type $(r,n)$ est de dimension $r+1$. C'est une cons\'equence directe, 
compte tenu  de sa d\'efinition, de la Proposition \ref{fond}. 

La m\'ethode canonique permet alors de donner
une repr\'esentation canonique d'un tissu semi-extr\'emal $\cal{T}$
comme tissu d'incidence sur un germe de vari\'et\'e, d'\'el\'ements
des courbes rationnelles normales de la vari\'et\'e de Blaschke de $\cal{T}$.
Mais nous ne savons pas si cette repr\'esentation est alg\'ebrique 
en g\'en\'eral. Toutefois dans le cas d'un tissu qui de plus est de rang
maximal, nous v\'erifions que sa vari\'et\'e de Blaschke 
appartient \`a une classe $\cal{X}_{r+1,n}(q)$, voir 
la D\'efinition \ref{r-n-q}. 

Ceci nous permet de traiter le cas 
des tissus de type $(r,2)$ sans difficult\'e dans la Section 3,
car une vari\'et\'e de la classe $\cal{X}_{r+1,2}(q)$ est une vari\'et\'e
de Veronese d'ordre $q$.

Les vari\'et\'es de Veronese jouent un r\^ole majeur 
dans \cite{P-Tr} o\`u elles apparaissent comme des projections 
distingu\'ees 
de vari\'et\'es d'une classe $\cal{X}_{r+1,n}(q)$
g\'en\'erale. Les plus simples des r\'esultats de Pirio-Tr\'epreau \cite{P-Tr},
expos\'es et compl\'et\'es dans les Sections 4 et 5,
suffisent alors pour montrer que le mod\`ele canonique d'un tissu de rang 
maximal est alg\'ebrique.

Nos autres \'enonc\'es principaux d\'ecoulent directement du 
pr\'ec\'edent, mais en s'appuyant sur des r\'esultats plus profonds
de \cite{P-Tr}.

\ssct{Propri\'et\'es de la vari\'et\'e de Blaschke}

Soit $\cal{T}$ un $d$-tissu semi-extr\'emal de type $(r,n)$
sur un germe $(M,x_0)$, de rang~$N+1$.
Rappelons, voir la D\'efinition~\ref{var-X},
que sa {\em vari\'et\'e de Blaschke} $X_{\cal{T}}$ 
est l'intersection de toutes les sous-vari\'et\'es projectives de $\P^N$
qui contiennent les courbes rationnelles normales $\s(x)$, $\,x\in (M,x_0)$,
de degr\'e $q(d)$. Nous notons $q=q(d)$.

\sk
Nous reprenons les notations de la Proposition \ref{fond}.
L'application de Blaschke, qui param\`etre 
la r\'eunion des courbes $\s(x)$
est de rang constant $r+1$. Son image $\tilde{X}$
est donc un germe, le long de la courbe $\s(x_0)$, 
de vari\'et\'e analytique 
lisse de dimension $r+1$. Soit 
$$
\k^{(n)} : (M,x_0)\times \left(\prod_{j=1}^n (\C,\t_j(x_0))\right)  \rightarrow \prod_{j=1}^n(\tilde{X},\k_j(x_0))
$$
l'application d\'efinie par $\k^{(n)}(x,t_1,\ldots,t_n) = (\k(x,t_1),\ldots,\k(x,t_n))$.

Le noyau de sa diff\'erentielle au point de base est donn\'e par les \'equations
$$
dt_j + J(x_0,\t_j(x_0))=0, \qquad \sum_{\a=1}^n \t_j(x_0)^{\a-1} m_{a,\a}(x_0)=0,
$$
o\`u $a=1,\ldots,r$ et  $j=1,\ldots,n$. Le second groupe d'\'equations donne $m_{a,\a}(x_0)=0$
pour tout $a$ et $\a$ donc $dx=0$. Le premier groupe d'\'equations
donne alors $dt_1=0, \ldots,dt_n=0$. Comme la source et le but ont la m\^eme dimension
$(r+1)n$, on obtient que $\k^{(n)}$ est un isomorphisme.

En particulier, pour tout voisinage $U$ dans $\tilde{X}$, par exemple du  
point $\k_1(x_0)$,
les courbes $\s(x)$ qui passent 
par ce point recouvrent un ouvert non vide de $U$.
\ble
La vari\'et\'e de Blaschke $X_\cal{T}$ du tissu $\cal{T}$
est une sous-vari\'et\'e projective irr\'eductible
de dimension $r+1$ de $\P^N$.
\ele
\bpf
Nous pouvons supposer que la remarque qui pr\'ec\`ede l'\'enonc\'e s'applique 
en $0\in \C^N$. Il s'agit de montrer que le germe lisse $(\tilde{X},0)$
est contenu dans une sous-vari\'et\'e alg\'ebrique de $\C^N$, de m\^eme dimension
$r+1$. L'irr\'eductibilit\'e de la vari\'et\'e $X_\cal{T}$ est alors une cons\'equence 
de sa d\'efinition et du fait que le germe $\tilde{X}$ est irr\'eductible.

La d\'emonstration qui suit est standard. (C'est l'occasion de r\'eparer 
une b\'evue dans \cite{Tr}, page 430.)
Nous param\'etrons l'espace des $(N+1)$-uplets 
de polyn\^omes $(p_1(t),\ldots,p_{N+1}(t))$,
nuls en $0$ et de degr\'e $\leq q$  par la famille $\g \in \C^{(N+1)q}$
de leurs coefficients. \`A tout $\g$ est associ\'e le germe d'application 
$$
t\in (\C,0), \qquad x_\g(t) = \left(\fr{p_1(t)}{1+p_{N+1}(t)}, \ldots, \fr{p_N(t)}{1+p_{N+1}(t)}\right).
$$
Quand $\g$ d\'ecrit $\C^{(N+1)q}$, on obtient en particulier 
des param\'etrages
de tous les germes lisses en $0$ de courbes rationnelles de degr\'e
$\leq q$.

\sk
Comme $p_{N+1}(0)=0$, les composantes de $x_\g(t)$ sont d\'eveloppables 
en des s\'eries enti\`eres de $t$ dont les coefficients  
sont des polyn\^omes en $\g$.

\sk
Si $f(x) = \sum_{|\a|\geq 1} a_\a x^\a$ est une fonction analytique 
sur $(\C^N,0)$ fix\'ee, comme $x_\g$ est de valuation $\geq 1$,
les coefficients de la s\'erie enti\`ere compos\'ee $f(x_\g(t))$
sont encore des polyn\^omes en $\g$. Cette s\'erie est donc nulle
si et seulement si $\g$ appartient \`a un ensemble alg\'ebrique 
d\'etermin\'e par la fonction $f(x)$.

\sk
Finalement, comme le germe $(\tilde{X},0)$ est d\'efini par l'annulation d'une famille 
de telles fonctions $f(x)$, on obtient qu'il existe un ensemble alg\'ebrique 
$K\subset \C^{(N+1)q}$ tel que l'image d'un germe 
$x_\g$ est contenue dans $(\tilde{X},0)$
si et seulement si $\g$ appartient \`a $K$.

\sk
Soit $\wh{K}$ la sous-vari\'et\'e projective de $\P^{(N+1)q}$
obtenue en compl\'etant $K$. Le germe d'application $(\g,t) \mapsto x_\g(t)$
est la restriction d'une application rationnelle 
$\phi: \; \wh{K}\times \P^1 \dasharrow \P^N$, 
dont l'image est une sous-vari\'et\'e 
projective $X$ de $\P^N$ qui contient par hypoth\`ese le germe de $\tilde{X}$
en $0$ et donc, par analyticit\'e, le germe $\tilde{X}$ lui-m\^eme.
Quitte \`a remplacer $\wh{K}$ par l'une de ses composantes irr\'eductibles,
bien choisie, nous pouvons, en conservant la propri\'et\'e pr\'ec\'edente,
supposer que $\wh{K}$ et $X$ sont irr\'eductibles.

\sk
Il reste \`a voir que $X$, qui est de dimension $\geq r+1$
par hypoth\`ese, est de dimension $\leq r+1$, autrement dit 
que $\phi$ est de rang $\leq r+1$ au point g\'en\'erique 
de $\wh{K}\times \P^1$. C'est \'evident car, si $U$ est un voisinage  assez petit 
d'un point donn\'e de $K$ et si $\epsilon > 0$ est assez 
petit, l'image de $U\times \{t\in \C, \; |t|<\epsilon\}$
par $\phi$ est contenue dans dans un voisinage donn\'e de $0$ dans $(\tilde{X},0)$,
par d\'efinition m\^eme de $K$.
\epf
L'ensemble des 
courbes rationnelles normales de degr\'e $q$ de $\P^N$ est un ouvert irr\'eductible
de la vari\'et\'e de Chow des $1$-cycles effectifs de degr\'e $q$ 
de $\P^N$. Notons $\text{CR}_{q}(\P^N)$ son adh\'erence 
et $\text{CR}_{q}(X_\cal{T})$ la vari\'et\'e alg\'ebrique ferm\'ee
des \'el\'ements de $\text{CR}_{q}(\P^N)$ dont le support est contenu
dans $X_\cal{T}$.
La vari\'et\'e d'incidence 
$$
I^{(n)} = \{(x,(a_1,\ldots,a_n))\in \text{CR}_{q}(X_\cal{T})\times X_\cal{T}^n, \; (a_1,\ldots,a_n)\in x^n\}
$$
est une vari\'et\'e projective compacte et la propri\'et\'e de 
l'application $\k^{(n)}$
mentionn\'ee au d\'ebut de cette section  
montre que la projection canonique $I^{(n)}\rightarrow X_\cal{T}^n$ est une submersion 
en au moins un point de $I^{(n)}$. Elle est donc surjective.

\sk
Ceci d\'emontre le Th\'eor\`eme \ref{Th-X0}. Nous le pr\'ecisons maintenant 
dans le cas o\`u le tissu $\cal{T}$ est de rang maximal.
\ble
La vari\'et\'e de Blaschke d'un $d$-tissu de rang maximal, de type $(r,n)$ 
avec $q(d)\geq n-1$, appartient \`a la classe $\cal{X}_{r+1,n}(q(d))$.
\ele
\bpf
C'est une cons\'equence directe de la D\'efinition~\ref{r-n-q}.
Il s'agit en fait de v\'erifier ici que cette d\'efinition co\"{\i}ncide 
avec celle qui est donn\'ee par la D\'efinition~1.1 de \cite{P-Tr}.
Les deux d\'efinitions diff\`erent seulement par la formule qui donne
la dimension $N$ de l'espace engendr\'e par une vari\'et\'e $X$ 
de la classe $\cal{X}_{r+1,n}(q)$.

Nous avons d\'efini ici $N$ par $N+1=\rho_{r,n}(d)$ o\`u $d=q+r(n-1)+2$.
Pour \'economiser l'espace, nous notons ${\text{C}}_{k}^{l}$ au lieu de $k \choose l$ les coefficients 
du bin\^ome. Soit 
$$
q = \rho(n-1) + m - 1, \qquad \rho\geq 1, \qquad m\in \{1,\ldots,n-1\},
$$
la division euclidienne de $q$ par $n-1$. Nous avons 
\beqn
N+1 
&  = & 
\sum_{h=0}^\rho \, 
((\rho-h)(n-1) + m) {\text{\large C}}_{r+h-1}^{r-1} \\
&  = & 
(\rho(n-1)+m) \sum_{h=0}^\rho \, {\text{\large C}}_{ r+h-1}^{r-1} 
- (n-1)\sum_{h=0}^\rho \, h{\text{\large C}}_{r+h-1}^{r-1},
\eeqn
et comme $\,h{\text{\large C}}_{r+h-1}^{r-1} = r {\text{\large C}}_{r+h-1}^{r}$
et $\,\sum_{h=0}^\rho {\text{\large C}}_{r+h-1}^{r-1} = {\text{\large C}}_{r+\rho}^{r}$,
\beqn
N+1
& = & (\rho(n-1)+m) {\text{\large C}}_{r+\rho}^{r} - r(n-1){\text{\large C}}_{r+\rho}^{r+1} \\
& = &
(r+1)(n-1) {\text{\large C}}_{r+\rho}^{r+1} + m{\text{\large C}}_{r+\rho}^{r} - r(n-1){\text{\large C}}_{r+\rho}^{r+1} \\
& = & m{\text{\large C}}_{r+\rho+1}^{r+1}  + (n-1-m){\text{\large C}}_{r+\rho}^{r}.
\eeqn
C'est la formule qui donne $N$ dans \cite{P-Tr}~Th\'eor\`eme~1.1.
\epf

\ssct{Le tissu $\cal{T}$ comme tissu d'incidence}

Nous continuons avec les notations et les hypoth\`eses de la section pr\'ec\'edente.
Nous rappelons d'abord ce que nous a donn\'e la m\'ethode canonique 
dans le Chapitre 3. Soit $\cal{F}_1,\ldots,\cal{F}_d$ les feuilletages qui composent le tissu 
semi-extr\'emal $\cal{T}$. Nous notons encore $q$ pour $q(d)=d-r(n-1)-2$.

\sk
Soit $(\phi_1^{(\l)},\ldots,\phi_d^{(\l)})$,
$\, 1\leq \l \leq N+1$, une base de l'espace des relations ab\'eliennes 
de $\cal{T}$ en $x_0$. Les applications de Poincar\'e sont donn\'ees par 
$$
x\in (M,x_0), \qquad {\k}_j (x) = [\phi_j^{(1)}(x): \cdots : \phi_j^{(N+1)}(x)].
$$
Chacune d'elles est de rang constant $r$ et 
constante le long des feuilles de $\cal{F}_j$.
Son image est donc un germe de vari\'et\'e lisse $(Z_j,p_j)$,
en $p_j=\k_j(x_0)$,  
de dimension $r$, et la submersion induite
$$
\k_j : (M,x_0) \rightarrow (Z_j, p_j)
$$
d\'efinit le feuilletage $\cal{F}_j$ : ses fibres sont les feuilles de $\cal{F}_j$.
Les points $p_j$ sont deux-\`a-deux distincts.

\sk
D'autre part, chaque point $x\in (M,x_0)$ d\'etermine une  
courbe $\s(x)$, une courbe rationnelle normale de degr\'e $q$, qui passe par les 
$d$ points $\k_j(x)$ et l'application $x\mapsto \s(x)$, \`a valeurs
dans une vari\'et\'e de Chow convenable, est une immersion,
donc induit un isomorphisme
$$
\s: (M,x_0) \rightarrow (\Si',\s(x_0)),
$$
o\`u  $(\Si',\s(x_0))$ est un germe de vari\'et\'e lisse de dimension $rn$.
Enfin, nous venons de montrer que 
ces courbes sont contenues dans une vari\'et\'e projective $X_\cal{T}$
de dimension $r+1$ dont le Th\'eor\`eme \ref{Th-X0} pr\'ecise les propri\'et\'es.

\bk
Nous identifions \`a pr\'esent $(M,x_0)$ \`a $(\Si',\s(x_0))$ par 
l'isomorphisme $\s$ et  nous consid\'erons le tissu
$\cal{T}=\{\cal{F}_1,\ldots,\cal{F}_d\}$ comme un tissu semi-extr\'emal de type 
$(r,n)$ sur $(\Si',x_0)$. 

Pour plus de clart\'e nous notons 
maintenant $x$ un point de $(\Si',x_0)$ et $\s(x)$ la courbe de $X_\cal{T}$
qu'il d\'efinit. 

La courbe $\s(x_0)$ est transverse aux germes $(Z_j,p_j)$ et le feuilletage 
$\cal{F}_j$ est d\'efini par la submersion 
d'incidence $\k_j: (\Si',x_0) \rightarrow (Z_j,p_j)$.

\sk
Nous avons donc la repr\'esentation suivante, 
canonique \`a un automophisme de $\P^N$
pr\`es, d'un tissu semi-extr\'emal de rang $N+1$.
\bpr
\label{rep1}
Soit $X$ la vari\'et\'e de Blaschke d'un $d$-tissu semi-extr\'emal $\cal{T}$,
de type $(r,n)$ avec $r\geq 2$ et de rang $N+1$.
Il existe  un germe $(\Si',x_0)$ de vari\'et\'e analytique lisse 
de dimension $rn$, de courbes rationnelles normales de degr\'e $q=q(d)$
de $X$, tel que $\cal{T}$ est isomorphe \`a {\em un tissu d'incidence}
sur $(\Si',x_0)$ pour la paire $(X,\Si')$.
\epr
On entend par l\`a que le tissu sur $(\Si',x_0)$ est d\'efini 
par les propri\'et\'es d'incidence des courbes  $\s(x)$ 
avec $d$ germes d'hypersurfaces lisses $(Z_j,p_j)$
de $X$, transverses \`a $\s(x_0)$ en des points $p_j$
deux-\`a-deux distincts. Nous ne disons pas qu'une telle 
configuration d\'efinit toujours un tissu, mais c'est le cas
par construction dans le cadre de l'\'enonc\'e.

\sk
Ce r\'esultat n'est pas satisfaisant. On {\em aimerait}
montrer que le germe $(\Si',x_0)$ est contenu dans une vari\'et\'e
alg\'ebrique de m\^eme dimension, ce que nous ne savons pas faire,
puis que les germes $(Z_j,p_j)$
sont contenus dans une hypersurface alg\'ebrique de $X_\cal{T}$.

\bk
Dans la suite de ce chapitre, nous ne consid\'erons plus que 
des tissus {\em de rang maximal}.

\sk
Soit $X$ une vari\'et\'e de la classe $\cal{X}_{r+1,n}(q)$.
On  montre, voir \cite{P-Tr}~Lem\-me~2.8\footnote{Nous n'utilisons pas
les m\^emes notations que dans \cite{P-Tr} o\`u cette vari\'et\'e 
est not\'ee $\ov{\Si}_q(X)$ et o\`u ce que nous noterons $\Si_{X,{\rm adm}}$
est not\'e $\Si_q(X)$.},
qu'il existe une et une seule composante irr\'eductible $\Si_X$
de l'adh\'erence de Zariski (dans une vari\'et\'e de Chow convenablement choisie)
de la vari\'et\'e des courbes rationnelles de degr\'e $q$  contenues 
dans $X$ telle que, 
pour tout $(a_1,\ldots,a_n)\in X^n$, il existe  un \'el\'ement 
de cette composante dont le support dans $X$ contient les points $a_j$.
Elle est de dimension $rn$.

\sk
Ceci permet de pr\'eciser la Proposition \ref{rep1} lorsque le tissu 
$\cal{T}$ est de rang maximal.
\bpr
\label{rep2}
Sous les hypoth\`ese de la Proposition \ref{rep1} et si le tissu $\cal{T}$
est de rang maximal, le germe $(\Si',x_0)$
est contenu dans la vari\'et\'e alg\'ebrique $\Si_X$,
de m\^eme dimension $rn$.
\epr

\ssct{Le cas des tissus de type $(r,2)$}

Nous d\'emontrons ici le Corollaire \ref{Th1D}, c'est-\`a-dire le Th\'eor\`eme \ref{Th1A}
pour $n=2$.

C'est un cas particulier remarquable. La d\'emonstration n'utilise rien de \cite{P-Tr}
\`a l'exception de la remarque suivante, qui 
peut remplacer un argument fautif \`a la fin de la d\'emonstration 
de H\'enaut \cite{He}.

Si $q\geq 1$, par d\'efinition, une vari\'et\'e $X$ de la classe $\cal{X}_{r+1,2}(q)$
engendre un espace de dimension $N={r+1+q\choose r+1}-1$ et a la propri\'et\'e
qu'une paire g\'en\'erale de points de $X$ est contenue dans une courbe rationnelle 
normale de degr\'e $q$ contenue dans $X$.

Ces propri\'et\'es caract\'erisent la {\em  vari\'et\'e de Veronese} de dimension $r+1$
et d'ordre $q$, c'est-\`a-dire l'image de $\P^{r+1}$ par le plongement de 
Veronese, d\'efini par le syst\`eme lin\'eaire $|\cal{O}_{\P^{r+1}}(q)|$.
Cette caract\'erisation appara\^{\i}t d\'ej\`a 
dans Bompiani \cite{Bo} avec une  d\'emonstration sans 
doute insuffisante. Elle est maintenant
bien \'etablie, voir le Th\'eor\`eme~1.5 de \cite{P-Tr}
et les commentaires qui le suivent.

\sk
La Proposition \ref{rep2}  nous permet de  
supposer que le $d$-tissu $\cal{T}$, de type $(r,2)$
et de rang maximal avec $q=q(d)\in \N^\star$, est un tissu d'incidence
pour la paire $(X_\cal{T},\Si_{X_\cal{T}})$.
Il poss\`ede en particulier une relation ab\'elienne compl\`ete,
voir plus bas la d\'emonstration du Th\'eor\`eme \ref{Th1B}.

Comme l'isomorphisme de Veronese 
$v_q : \P^{r+1} \rightarrow X_\cal{T}$ envoie 
isomorphiquement la paire 
$(\P^{r+1},\G_{r,2})$ sur la paire $(X_\cal{T},\Si_{X_\cal{T}})$,
le tissu est isomorphe \`a un tissu d'incidence 
pour la paire $(\P^{r+1},\G_{r,2})$.
Compte tenu du th\'eor\`eme d'Abel inverse, celui-ci est alg\'ebrique grassmannien.

\ssct{Tissus d'incidence pour une paire $(X,\Si_X)$}

Soit $X$ une vari\'et\'e de la classe $\cal{X}_{r+1,n}(q)$.
Nous pr\'ecisons maintenant la construction de tissus 
de type $(r,n)$ esquiss\'ee dans la Section 1.4.
Concernant les propri\'et\'es de la vari\'et\'e $X$,
nous nous r\'ef\'erons au Chapitre~2 de \cite{P-Tr}.

\sk
Soit 
$$
q=\rho(n-1) + m-1, \qquad \rho\geq 1, \qquad m\in\{1,\ldots,n-1\},
$$
la division euclidienne de $q$ par $n-1$.

\sk
Notons $X_a(k)$ le sous-espace de $\P^N=\lan X\ran$, osculateur \`a l'ordre $k$
\`a $X$ en $a\in X_{\rm reg}$. 
Soit $(a_1,\ldots,a_{n-1})\in X_{\rm reg}^{n-1}$ et 
$$
E_i = X_{a_i}(\rho), \;\; 1\leq i\leq m-2 \; ; \;\; E_i=X_{a_i}(\rho-1), \;\; m-1\leq i\leq n-1.
$$
On montre que pour $(a_1,\ldots,a_{n-1})$ g\'en\'erique, on a 
la d\'ecomposition 
\beq
\label{directe}
\P^N = \oplus_{i=1}^{n-1} E_i
\eeq
de $\P^N$ en somme directe projective des $E_i$ :
ces espaces engendrent $\P^N$ et $\sum_{i=1}^{n-1} (\dim E_i + 1)=N+1$.

On lui associe la {\em projection osculatrice} de centre $\oplus_{i=2}^{n-1}E_i$ :
\beq
\label{proj}
\t : \P^N\bck \oplus_{i=2}^{n-1} E_i \rightarrow E_1.
\eeq
Bien s\^ur, si $n=2$, le centre de la projection est vide et celle-ci est l'identit\'e.

\sk
Nous avons introduit dans \cite{P-Tr} la propri\'et\'e d'\^etre 
admissibles pour diff\'erentes sortes d'objets.
La notion de base est celle d'un {\em $n$-uplet admissible}
$(a_1,\ldots,a_n)\in X^n$. C'est une famille de $n$ points 
deux-\`a-deux distincts de $X_{\rm reg}$ telle 
que $n-1$ quelconques d'entre eux, pris dans 
n'importe quel ordre, d\'efinissent
une d\'ecomposition de $\P^N$ analogue \`a (\ref{directe}).

\sk
Tout $n$-uplet admissible est contenu dans une 
courbe rationnelle normale de degr\'e $q$
d\'efinie par un point $x$ de $\Si_X$. Une telle 
courbe est appel\'ee une {\em courbe admissible} de $X$. 

\bk
La projection $\t: \P^N \dasharrow X_{a_1}(\rho)$,
d\'efinie \`a partir de $n-1$ points $a_1, a_2,\ldots,a_{n-1}$
d'un $n$-uplet admissible, contenu dans une courbe admissible $\s(x_0)$, 
a les propri\'et\'es remarquables suivantes\footnote{Voir \cite{P-Tr} Chapitre 2.
Les objets admissibles sont introduits dans les D\'efinitions 2.3 et 2.4.
Les premi\`eres propri\'et\'es de la projection $\t$
apparaissent dans la Proposition~2.6. La troisi\`eme appara\^{\i}t 
plus tard, dans le Th\'eor\`eme 2.11.},
que nous mentionnons en premier car elles nous servirons 
pour d\'emontrer un th\'eor\`eme d'Abel inverse \`a la fin de cette section.

\be
\item  La projection $\t$ induit une application birationnelle de $X$
sur son image~$X'$, qui est une vari\'et\'e de la classe 
$\cal{X}_{r+1,2}(\rho)$, donc  {\em une vari\'et\'e 
de Veronese d'ordre $\rho$}.

\sk
\item  Pour tout $a\in \s(x_0)\bck\{a_2,\ldots,a_{n-1}\}$,
la projection $\t$ induit un isomorphisme du germe $(X,a)$ sur 
le germe $(X',\t(a))$.

\sk
\item Si $x\in (\Si_X(a_2)\cap \cdots \cap \Si_X(a_{n-1}),x_0)$, 
son image par $\t$ est une courbe rationnelle normale de degr\'e $\rho$ de $X'$
et l'application ainsi d\'efinie 
$$
\tau: (\Si_X(a_2)\cap \cdots \cap \Si_X(a_{n-1}),x_0) \rightarrow (\Si_{X'},\tau(x_0)),
$$
est un isomorphisme.  
\ee
Nous avons not\'e $\Si_X(a)$ l'ensemble alg\'ebrique des $x\in \Si_X$
dont le support contient le point $a$ de $X$.

Un isomorphisme de Veronese 
envoie $(\Si_{X'},\tau(x_0))$ sur un germe de la grassmannienne des droites 
de $\P^{r+1}$. On en d\'eduit que les courbes $\s(x)$, $x$ voisin de $x_0$,
qui passent par $a_2,\ldots,a_{n-1}$, recouvrent $(X,a_1)$
et que les droites tangentes $T_{a_1}\s(x)$
de celles qui passent aussi par le point $a_1$, 
recouvrent  un voisinage de $T_{a_1}\s(x_0)$ dans $T_{a_1}X$.

\sk
Les  propri\'et\'es qui suivent sont 
des cons\'equences 
des pr\'ec\'edentes, voir \cite{P-Tr} Sections 2.5 et 2.6.
L'ensemble des $x\in \Si_X$ tels que la courbe 
$\s(x)$ est admissible est un ouvert lisse et Zariski-dense 
$\Si_{X,{\rm adm}}$ de $\Si_X$. Tout $n$-uplet de points distincts 
d'une courbe admissible est admissible. L'ensemble 
des points (dits admissibles) de $X$ qui 
appartiennent \`a au moins un $n$-uplet admissible,
est un ouvert lisse et Zariski-dense  $X_{\rm adm}$
de $X$. C'est la r\'eunion de toutes les courbes 
admissibles de $X$.

Enfin, si $a_1,\ldots,a_n$ sont $n$ points distincts
d'une courbe admissible $\s(x_0)$,
les germes $(\Si_X(a_j),x_0)$ sont lisses 
de codimension $r$ et se coupent proprement : 
le germe en $x_0$ de $\Si_X(a_1)\cap \dots \cap \Si_X(a_p)$ est lisse de codimension $pr$,
d'espace tangent en $x_0$ l'intersection $T_{x_0}\Si_X(a_1)\cap\cdots \cap T_{x_0}\Si_X(a_p)$.

\sk
Soit $\s(x_0)$ une courbe admissible et $(Z,p)$ 
un germe d'hypersurface lisse transverse \`a $\s(x_0)$
en $p$. Compte tenu des propri\'et\'es
qu'on vient de rappeler, 
le morphisme d'incidence $\k : (\Si_X,x_0)\rightarrow (Z,p)$
est une submersion et donc d\'efinit un feuilletage
de type $(r,n)$ sur $(\Si_X,x_0)$, dont les feuilles 
sont des morceaux de vari\'et\'es de la forme 
$\Si_X(a)$.

Remarquons que si le germe $(Z,p)$
\'etait singulier ou s'il n'\'etait pas transverse 
\`a $\s(x_0)$, une courbe admissible voisine de $\s(x_0)$
g\'en\'erale rencontrerait le germe en plusieurs points et l'incidence ne 
d\'efinirait pas un feuilletage (r\'egulier).

R\'eciproquement, si $\cal{F}$ est un feuilletage 
de type 
$(r,n)$ sur $(\Si_X,x_0)$ dont les feuilles sont des morceaux 
de vari\'et\'es de la forme $\Si_X(a)$,
ces feuilles d\'efinissent un germe 
d'hypersurface $(Z,p)$ avec $x_0=\Si_X(p)$ et compte tenu de ce 
qu'on vient de dire, ce germe est lisse et transverse \`a $\s(x_0)$.

\sk
Nous avons donc deux d\'efinitions \'equivalentes de ce que 
nous appelons un {\em feuilletage d'incidence admissible},
en un point de $\Si_{X,{\rm adm}}$.

\sk
Si $x_0\in \Si_{X,{\rm adm}}$, on obtient un {\em $d$-tissu admissible} sur $(\Si_X,x_0)$
en se donnant $d$ germes lisses $(Z_j,p_j)$
d'hypersurfaces transverses \`a la courbe $\s(x_0)$
en des points distincts.
Il s'agit bien d'un tissu, puisque $p\leq n$ parmi les vari\'et\'es 
$\Si_X(p_j)$ se coupent proprement en $x_0$, comme on a dit.
(Voir la Proposition \ref{rep1}
pour une notion peut-\^etre plus g\'en\'erale 
de tissu d'incidence.)

\sk
Nous pouvons donner un contenu rigoureux \`a la D\'efinition \ref{Conv2}.
\ble
\label{tissu-alg}
Soit $Z$ une hypersurface alg\'ebrique r\'eduite de $X$, 
dont toutes les composantes 
irr\'eductibles rencontrent $X_{\rm adm}$.
Pour $x_0\in \Si_{X,{\rm adm}}$, elle d\'efinit
un tissu admissible sur le germe $(\Si_X,x_0)$,
d'ordre son nombre d'intersection avec toute courbe 
admissible de $X$.
\ele 
Nous dirons que ces tissus sont des germes du {\em tissu alg\'ebrique 
d'incidence} $\cal{T}_Z$ d\'efini par l'hypersurface $Z$.
\bpf
Il suffit de v\'erifier qu'une courbe admissible g\'en\'erale 
rencontre $Z$ transversalement en $d$ points distincts 
de $Z_{\rm reg}$. C'est une cons\'equence des Propri\'et\'es 1, 2, 3 et
du commentaire qui suit leur \'enonc\'e : il montre qu'en partant 
d'une courbe admissible quelconque, des perturbations 
successives de cette courbe permettent de perturber
ses points d'intersection avec $Z$ ainsi que 
les directions de la courbe en ces points.
\epf

Nous avons la version suivante du th\'eor\`eme d'Abel inverse 
pour la paire $(X,\Si_X)$.
\bt
\label{4-abel}
Soit $\cal{T}$ un $d$-tissu d'incidence admissible pour la paire $(X,\Si_X)$.
S'il poss\`ede une relation ab\'elienne compl\`ete,
c'est un germe d'un tissu alg\'ebrique d'incidence pour cette paire,
d\'efini par une hypersurface alg\'ebrique $Z$ de $X$.
Ses relations ab\'eliennes sont induites par les  $r$-formes $\Si_X$-ab\'eliennes sur $Z$
et celles-ci sont rationnelles.
\et
\bpf
Notons $\k_j: (\Si_X,x_0)\rightarrow (Z_j,p_j)$ les morphismes d'incidence qui 
d\'efinissent le tissu $\cal{T}$ de l'\'enonc\'e.
Par hypoth\`ese, il existe 
des $r$-formes lisses $\phi_j$ sur les germes $(Z_j,p_j)$, non nulles et
dont la trace $\sum_{j=1}^d \k_j^\star \phi_j$ est la forme 
nulle.

\sk
Si $n=2$, le th\'eor\`eme se r\'eduit au th\'eor\`eme d'Abel inverse usuel
puisque la paire $(X,\Si_X)$ est alors isomorphe \`a la paire $(\P^{r+1},\G_{r,2})$.
Nous supposons maintenant $n\geq 3$. 

\sk
Nous nous ramenons au cas $n=2$ 
gr\^ace \`a une projection osculatrice $\t: X\dasharrow X'$
d\'efinie par des points $a_1,\ldots,a_{n-1}$ de 
$\s(x_0)$, autres que les $p_j$.
Nous utilisons les Propri\'et\'es 1, 2 et 3 ci-dessus.
La vari\'et\'e $X'$ est une vari\'et\'e de Veronese d'ordre $\rho$.
Consid\'erons la situation obtenue sur $X'$
par projection des donn\'ees. 

\sk
Sur $X'$ nous avons $d$ germes d'hypersurfaces lisses  $(Z'_j,p'_j)$,
transverses \`a la courbe $x'_0=\tau(x_0)$ aux points deux-\`a-deux 
distincts $p'_j$,  et sur chaque $(Z_j,p_j)$ 
la $r$-forme lisse non nulle $\phi'_j= ({\t_{|Z_j}})_\star \phi_j$.

\sk
Notons 
$$
\k'_j : (\Si_{X'},x'_0) \rightarrow (Z'_j,p'_j)
$$
les morphismes d'incidence correspondant \`a la nouvelle situation.

\sk
Dans la situation initiale, les morphismes d'incidence $\k_j$ induisent
aussi des morphismes d'incidence 
$$
\l_j: (\Si_X(a_1)\cap \cdots \cap \Si_X(a_{n-2}),x_0) \rightarrow (Z_j,p_j),
$$
o\`u $\l_j = \k_j\circ \iota$ et $\iota$ est l'inclusion de 
$\Si_X(a_1)\cap\cdots\cap\Si_X(a_{n-2})$ dans $\Si_X$.
Par hypoth\`ese, la trace $\sum_{j=1}^d \l_j^\star \phi_j = \iota^\star (\sum_{j=1}^d \k_j^\star \phi_j)$
est nulle. 

Comme $\k'_j = \t_{|Z_j}\circ  \l_j \circ \tau^{-1}$ pour tout $j$, 
nous obtenons que la trace 
$$
\sum_{j=1}^d {\k'_j}^\star \phi'_j 
=
(\tau^{-1})^\star \sum_{j=1}^d \l_j^\star \phi_j
$$
est nulle au voisinage de $\tau(x_0)$ dans $\Si_{X'}$.

\sk
Le th\'eor\`eme est vrai pour la vari\'et\'e de Veronese 
$X'$ donc les germes $(Z'_j,\tau(x_0))$ sont contenus dans 
une hypersurface alg\'ebrique $Z'$, que nous pouvons supposer 
minimale pour l'inclusion, et les formes $\phi'_j$ 
sont induites par une $r$-forme 
rationnelle $\phi'$ sur $Z'$.

\sk
Remontant \`a $X$, nous obtenons que l'image (stricte) 
de $Z'$ par $\t^{-1}$ est une hypersurface alg\'ebrique 
$Z$ de $X$ et $\phi = (\t_{|Z})^\star \phi'$ une $r$-forme 
rationnelle sur $Z$, telles que les germes $(Z_j,p_j)$
sont contenus dans $Z$ et que les formes $\phi_j$
sont induites par restriction de la forme $\phi$. 
\epf
\bco
Soit $Z$ une hypersurface r\'eduite de $X$, sans composante
irr\'eductible contenue dans $X\bck X_{\rm adm}$.
Les $r$-formes $\Si_X$-ab\'eliennes sur $Z$ sont rationnelles.
\eco
\'Etant donn\'e une $r$-forme $\Si_X$-ab\'elienne $\phi$ sur $Z$,
il suffit d'appliquer le r\'esultat pr\'ec\'edent au tissu 
alg\'ebrique $\cal{T}_{Z'}$, o\`u $Z'$ est la r\'eunion 
des composantes irr\'eductibles de $Z$ sur lesquelles 
$\phi$ n'est pas la forme nulle.

\ssct{D\'emonstration des principaux \'enonc\'es}

Nous achevons d'abord la d\'emonstration du Th\'eor\`eme \ref{Th1B}.
Les deux autres \'enonc\'es sont des cons\'equences 
de celui-ci et de \cite{P-Tr}. Le lemme suivant se d\'emontre comme le Lemme \ref{tissu-alg}.
Il montre que le th\'eor\`eme d'Abel inverse pr\'ec\'edent peut-\^etre appliqu\'e 
\`a un tissu d'incidence sur un germe $(\Si',x_0)$, voir la 
Proposition \ref{rep1}, m\^eme si $\s(x_0)$ n'est pas une 
courbe admissible de $X$.
\ble
Soit $(\Si',x_0)$ un germe lisse dans $\Si_X$,
de m\^eme dimension $rn$ que $\Si_X$, et $\cal{T}$ un tissu 
d'incidence pour la paire $(X,\Si')$.
Pour $x\in (\Si',x_0)$ g\'en\'eral, le tissu induit par $\cal{T}$
sur $(\Si',x)$ est admissible.
\ele

\bpf[D\'emonstration du Th\'eor\`eme \ref{Th1B}]

Nous pouvons supposer que le $d$-tissu  $\cal{T}$, de type $(r,n)$ et de rang maximal, 
de vari\'et\'e de Blaschke $X$, est d\'efini sur un germe lisse $(\Si',x_0)$ 
contenu dans $\Si_X$ et de m\^eme dimension.

Le tissu $\cal{T}$ admet une relation ab\'elienne {\em compl\`ete}.
En effet, il admet $q(d)+1$ relations ab\'eliennes 
dont les $0$-jets sont lin\'eairement ind\'ependants quand ses sous-tissus 
propres en admettent au plus $q(d)$. Une combinaison lin\'eaire 
convenable de ces relations est donc compl\`ete. 

Compte tenu du
Th\'eor\`eme~\ref{4-abel} et du lemme pr\'ec\'edent, le tissu $\cal{T}$ est un germe 
d'un tissu alg\'ebrique d'incidence pour la paire $(X,\Si_X)$,
ce qui donne le r\'esultat.
\epf

\bpf[D\'emonstration du Th\'eor\`eme \ref{Th1C}]

Nous continuons avec les hypoth\`eses et les notations de la d\'emonstration 
pr\'ec\'edente. 

Nous savons que le tissu $\cal{T}=\{\cal{F}_1,\ldots,\cal{F}_d\}$ d\'etermine 
une  ${\rm G}_{r,n}$-structure sur le germe lisse $(\Si',x_0)$ o\`u il est d\'efini,
dont les feuilles des feuilletages $\cal{F}_j$
sont des sous-vari\'et\'es distingu\'ees.
De plus, le tissu est isomorphe \`a un germe d'un tissu
alg\'ebrique grassmannien si et seulement si cette structure
est int\'egrable.  

D'autre part, voir \cite{P-Tr} Th\'eor\`eme 4.4, la vari\'et\'e $\Si_{X,{\rm adm}}$
admet une ${\rm G}_{r,n}$-structure, d\'etermin\'ee 
par le fait que les vari\'et\'es 
$$
\Si_X(p) = \{x\in \Si_X, \;\; p\in \s(x)\}, \qquad p\in X,
$$
induisent des sous-vari\'et\'es distingu\'ees sur $\Si_{X,{\rm adm}}$.
Il est clair, de par leurs caract\'erisations, 
que les deux structures ci-dessus co\"{\i}ncident 
sur $\Si'\cap \Si_{X,{\rm adm}}$, un ouvert dense de $(\Si',x_0)$.
Si la ${\rm G}_{r,n}$-structure de $\Si_{X,{\rm adm}}$ est int\'egrable,
la structure d\'efinie par le tissu $\cal{T}$ est int\'egrable
sur un ouvert dense de $(\Si',x_0)$ et donc au voisinage de $x_0$
car, comme  il s'agit d'une $\rm G$-structure de type fini, 
la propri\'et\'e d'\^etre int\'egrable au voisinage d'un point est
ferm\'ee. D'o\`u le th\'eor\`eme.
\epf

\bpf[D\'emonstration du Th\'eor\`eme \ref{Th1A}]
Le cas $n=2$ a d\'ej\`a \'et\'e trait\'e. 

Le th\'eor\`eme 
est alors une cons\'equence imm\'ediate d'un r\'esultat essentiel,
le  Corollaire 3.8 de \cite{P-Tr},
dont le sens est que si $X$ est une vari\'et\'e de la classe $\cal{X}_{r+1,n}(q)$
avec $q\neq 2n-3$, la ${\rm G}_{r,n}$-structure de $\Si_{X,{\rm adm}}$ 
est int\'egrable. 

Pour \^etre pr\'ecis, ce corollaire 
affirme que si $q\neq 2n-3$, la vari\'et\'e $X$ est \og int\'egrable \fg,
une notion provisoire. Elle est alors \og standard \fg\, d'apr\`es
\cite{P-Tr} Lemme 3.13, ce qui implique que la ${\rm G}_{r,n}$-structure de $\Si_{X,{\rm adm}}$
est int\'egrable d'apr\`es \cite{P-Tr}~Th\'eor\`eme 4.4.
\epf

\end{document}